\documentclass[11pt]{amsart}
\usepackage{amsthm}
\usepackage[bottom=2.8cm]{geometry}
\usepackage{amsmath,amssymb,mathrsfs}
\usepackage{graphicx}
\usepackage{enumerate}
\usepackage{mathtools}
\usepackage{xcolor}

\usepackage{dsfont}
\usepackage{hyperref}

\numberwithin{equation}{section}

\newtheorem{thm}{Theorem}[section]
\newtheorem{lem}[thm]{Lemma}
\newtheorem{prop}[thm]{Proposition}
\newtheorem{cor}[thm]{Corollary}

\newtheorem{notation}[thm]{Notation}

\newtheorem{rmk}[thm]{Remark}

\newcommand\supp{\mathrm{supp}\;}

\newcommand\halfopen[2]{\ensuremath{(#1,#2]}}

\newcommand{\N}{{\mathbb N}}
\newcommand{\R}{{\mathbb R}}

\newcommand{\C}{{\mathbb C}}

\author{Jacques Benatar}
\address{Brussels, Belgium }
\email{jbenatarmath@gmail.com}

\begin{document}

\begin{abstract} In the late eighties, Hildebrand and Tenenbaum proved an asymptotic formula for the number of positive integers below $x$, having exactly $\nu$ distinct prime divisors: $\pi_{\nu}(x) \sim x \delta_{\nu}(x)$.  Here we consider the restricted count $\pi_{\nu}(x,y)$ for integers lying in the short interval $\halfopen{x}{x+y}$. In this setting, we show that for any $\varepsilon >0$, the asymptotic equivalence
\[ \pi_{\nu}(x,y) \sim y \delta_{\nu}(x)\]
holds uniformly over all $1 \le \nu \le (\log x)^{1/3}/(\log \log x)^2$ and all $x^{17/30 + \varepsilon} \leq y \leq x$.  The methods also furnish mean upper bounds for the $k$-fold divisor function $\tau_k$ in short intervals, with strong uniformity in $k$. 
\end{abstract}

\title{A short-interval Hildebrand-Tenenbaum theorem }

\maketitle

\section{ Introduction and statement of the main results}

\subsection{Setting the stage}
The central problem of this article is to understand the short-interval density of the set $S_{\nu}$, which consists of those positive integers having exactly $\nu$ distinct prime factors.  To be precise, let $\omega(n)$ denote the number of distinct prime divisors of $n$ and define, for each positive integer $\nu $ and any $x\geq 1$, 
\[ S_{\nu}(x)= \left\{ n \leq x: \omega(n)=\nu \right\}, \qquad  \pi_{\nu}(x):=\left| S_{\nu}(x) \right|.  \]

Let us recall two classical results regarding the cardinality $ \pi_{\nu}(x)$. First we have the asymptotic formula of Landau \cite{Lan} which was obtained as early as 1900 and asserts that 
\begin{equation}\label{Lan}
 \pi_{\nu}(x) \sim \frac{x}{\log x} \frac{(\log_2 x)^{\nu -1}}{(\nu -1)!}
 \end{equation}
for any fixed $\nu$ (here we have used the notation $\log_2 x= \log \log x$). On the other hand, Hardy and Ramanujan \cite{HR} later gave a uniform upper bound, guaranteeing the existence of a constant $c>0$ such that
\begin{equation}\label{HRnu}
\pi_{\nu}(x) \ll \frac{x}{\log x} \frac{(\log_2 x+c)^{\nu -1}}{(\nu -1)!},  
\end{equation}
for any natural number $\nu$ and any $x\geq 3$. Refinements of the asymptotic \eqref{Lan} can be found in subsequent works of Sathe \cite{Sathe}, Selberg \cite{Sel}, Hensley \cite{Hen} and Pomerance \cite{Pom}, where the authors obtain sharper expansions for $\pi_{\nu}(x)$ as well as greater uniformity in $\nu$. 

There are a number of intriguing areas of inquiry connected with the study of $S_{\nu}$. Arguably, the two most prominent examples are the extensively studied \textsl{smooth numbers} (for general surveys on the topic as well as interactions with algorithmic number theory, we refer the reader to the article of Granville \cite{Gran} and the work of Hildebrand and Tenenbaum \cite{HTsurvey}) and the problem of approximating the cardinality $ \Pi_{\nu}(x)=|\left\{ n \leq x: \Omega(n)=\nu \right\}| $, where $\Omega$ counts prime divisors \textsl{with} multiplicity (see for instance the work of Hwang \cite{Hwang}). We also mention that Wolke and Zhan \cite{WZ} produced a Bombieri-Vinogradov type equidistribution result for arithmetic progressions in $S_{\nu}$.  

Returning to Landau's formula \eqref{Lan}, it is important to note that the asymptotic main term is no longer valid beyond the region $\nu \asymp  \log \log x$, instead taking on the appearance in \eqref{HilTen} just below.
The best result to date\footnote{In fact, Kerner obtained an improvement of the Hildebrand-Tenenbaum result in his thesis \cite{Kerner}.  The author gives an asymptotic formula for $\pi_{\nu}(x)$ in the very long range $\nu \leq (1+\frac{1-\varepsilon}{\log_2 x}) \frac{\log x}{\log_2 x}.$} is due to Hildebrand and Tenenbaum \cite{HT} and asserts that the natural density $\delta_{\nu}(x)$ of the set $S_{\nu}(x)$ is given by the asymptotic formula

\begin{equation}\label{HilTen}
\pi_{\nu}(x) = x \delta_{\nu}(x)  = \frac{ G(\rho,\alpha)x^{\alpha} \rho^{-\nu} }{ \nu \phi(\nu) \phi(\rho) \log x} \big( 1+ O(1/L_{\nu}(x))\big),
\end{equation}
provided that $1 \le \nu \le \log x/ (\log \log x)^2$. Here 
\begin{equation}\label{Gdef}
G(z,s)=\sum_{n \ge 1} z^{\omega(n)} n^{-s} =\prod_{p} \big(1+\frac{z}{p^s-1} \big) 
\end{equation}
and 
\begin{equation}\label{psiLdef}
 \phi(t)=\Gamma(t) t^{-t} e^{t}, \qquad L_{\nu}(x)=\log \Big( \frac{\log x}{\nu \log (\nu +1) } \Big),  
 \end{equation}
and the pair $(\rho,\alpha)$ is the unique minimiser of the quantity $G(r,a)x^{a} r^{-\nu}$ in the region $(r,a) \in (0,\infty) \times (1,\infty)$.

In the present work we seek to prove two short-interval results concerning $S_{\nu}$. Our first main goal  is to establish a short-interval version of \eqref{HilTen} for the restricted count
\begin{equation*}
\pi_{\nu}(x,y):=\left| S_{\nu}(x,y) \right|,  \qquad S_{\nu}(x,y)= S_{\nu} \cap \halfopen{x}{x+y}, 
\end{equation*}
with a strong emphasis on the uniformity in $\nu$.  The second purpose of this article is to shed some light on the anatomical structure of the \textsl{typical} integer $n \in S_{\nu}(x,y)$, which is to say, the size and distribution of its divisors.

\subsection{The anatomy of \boldmath\texorpdfstring{$S_{\nu}$}{S{nu}}.}  A useful \textsl{anatomical} feature of $S_{\nu}$ is that generic elements of the set are overwhelmingly made up of small prime factors. This phenomenon can be made precise by comparing the full count $\pi_{\nu}(x)$ to a suitable minorant. With this in mind, set $\lambda^{+}=  \log x / \sqrt{\log_3 x}$,  $\tau=\exp(\lambda^{+})$ and first consider the sums
\begin{equation}\label{restrictedlower}
 \mathcal{M}_{\nu}' (x)=\sum_{ m \in S_{\nu -1}(\tau) }\ 
\sum_{p  \leq x/m } 1, \qquad \ \  \mathcal{M}_{\nu}' (x,y)=\mathcal{M}_{\nu}' (x+y)-\mathcal{M}_{\nu}' (x).
\end{equation}

\noindent The minorant $ \mathcal{M}_{\nu}' (x)$ will be useful when $\nu=o(\log_2 x)$.  Indeed, in this regime the typical integer in $S_{\nu}(x)$ is expressible as the product of a single large prime divisor $p$ and a small divisor in $S_{\nu -1}(x/p)$.  

However, for larger values of $\nu$, the typical $n \in S_{\nu}(x)$ has a large divisor composed of roughly $\nu/L_{\nu}(x)$ many distinct primes.  In this general setting we thus require a variant of $\mathcal{M}_{\nu}' (x)$: set 
\begin{equation}\label{Mhashvariablesdef}
\ell_{\nu}(x) = \frac{ \nu \log(L_{\nu}(x))^2 }{L_{\nu}(x) }, \qquad \mathcal{Q}_{\nu}(x)=\prod_{p \leq t } p
\end{equation} 

with $t=\exp( \log x/ (\ell_{\nu}(x) \log_3 x) )$, and define the minorant
\begin{equation}\label{broadvariant}
 \mathcal{M}^{\sharp}_{\nu} (x)=\sum_{1 \leq w \leq \ell_{\nu}(x) } \ 
\sum_{ \substack{m \in S_{\nu -w}(\tau) \\ m | \mathcal{Q}_{\nu}(x)^{\infty}} }\ 
\sum_{\substack{ n  \in S_w( x/m) \\  (n, \mathcal{Q}_{\nu}(x) )=1 } } 1,
\end{equation}
together with its short-interval variant
\begin{equation}
 \mathcal{M}^{\sharp}_{\nu} (x,y)=  \mathcal{M}^{\sharp}_{\nu} (x+y) - \mathcal{M}^{\sharp}_{\nu} (x).
\end{equation}
The condition $m | \mathcal{Q}_{\nu}(x)^{\infty} $ appearing in the definition of $ \mathcal{M}^{\sharp}_{\nu} (x)$, is shorthand for the property $p|m \implies p \leq t $. 

Our first theorem,  stated  just below,  confirms the aforementioned anatomical structure of $S_{\nu}(x,y)$. In other words, almost all integers in $S_{\nu}(x,y)$ are captured by the minorant $ \mathcal{M}^{\sharp}_{\nu} (x,y)$. 

\begin{thm}\label{squarefreemainprop}
Let $\varepsilon>0$ and $a>4$ be fixed and let $\mathcal{L}_{a}(x)=\log x/(\log_2 x)^a $ for all sufficiently large $x \geq 1$.  Then for all $x^{17/30 + \varepsilon} \leq y \leq x$,  we have the estimate
\begin{equation}\label{shortintervalminorantupper}
  \pi_{\nu}(x,y) = (1+o(1))  \mathcal{M}^{\sharp}_{\nu}(x,y)  ,
\end{equation} 
uniformly in $\nu \leq \mathcal{L}_a(x)$.  Moreover,  in the shorter range $\nu \leq \log_2 x/ (\log_3 x)^{3/4}$, we have that
\begin{equation}\label{simpleanatomy}
  \pi_{\nu}(x,y) = (1+o(1))  \mathcal{M}^{'}_{\nu}(x,y).
\end{equation} 
\end{thm}

\subsection{Prior work on \texorpdfstring{$S_{\nu}(x,y)$}{S{nu}(x,y)}}\label{priorsection} For "small" values of $\nu$ we are fortunate to have a variety of striking short-interval results - some of these, such as Hoheisel's theorem on the existence of primes in short intervals, go back to the early twentieth century. A natural starting point for the study of $S_{\nu}(x,y)$ is the case $\nu=1$ which corresponds to the prime number theorem in short intervals. If one is interested in precise asymptotics of the form 
\begin{equation}\label{shortpnt}
 \psi(x,y)=\sum_{n \in [x,x+y]} \Lambda (n) \sim y,
 \end{equation}
the best result, until very recently, was due to Huxley \cite{Hux}, yielding an admissible range of $y\geq x^{7/12+ \varepsilon}$. This has been improved thanks to the work of Guth and Maynard \cite{GM} where the authors show that \eqref{shortpnt} remains valid for $y\geq x^{17/30 +\varepsilon}$. Further, in the case $\nu=2$, Matom\"{a}ki and Ter\"{a}v\"{a}inen \cite[Theorem 1.4]{MT}  recently proved the asymptotic $\pi_2(x,y)\sim y \delta_2(x)$ in the short range $y \geq x^{11/20+ \varepsilon}$.  

Alternatively, if one is satisfied with a lower bound  for $\pi_{\nu}(x,y)$, or if one is interested in almost-primes instead of $S_{\nu}$, sieve methods provide powerful tools to extend the range of $y$. The reader will find interesting theorems and references in the works of Baker-Harman-Pintz \cite{BHP} (for $\pi_1(x,y)$) and Wu \cite{Wu} (for $P_2$ almost primes). There is also a long history of results concerning \textsl{almost all} short intervals where, beginning with the work of Wolke \cite{Wolke}, various authors have studied very short intervals of the form $[x,x+(\log x)^c ]$ with $c>1$ (see the recent work of Ter\"{a}v\"{a}inen \cite{Ter} for the sharpest results in this direction).   

To conclude this introductory segment, we mention the case when $\nu$ is a slowly growing function of $x$. Here it is well known that, via the application of standard sieve methods, the formula \eqref{shortpnt} gives rise to an asymptotic expression for $\pi_{\nu}(x,y)$ (corresponding to Landau's formula \eqref{Lan}), provided that $\nu=o(\log_2 x)$.  Moreover, K\'atai \cite{Kat} and later Bassily-K\'atai \cite{BKat} obtained a formula for $\pi_{\nu}(x,y)$ in the range $\nu =O( \log_2 x)$ and $y \geq x^{7/12 + \varepsilon}$. In the preliminary sections below, we will briefly revisit some of the aforementioned results, before moving on to the unexplored range $\nu \gg \log_2 x$.

\subsection{Comparative results for \texorpdfstring{$\pi_{\nu}(x,y)$}{pi{nu}(x,y)} and estimates for short divisor sums}
Our second main result, just below, establishes comparative estimates/asymptotics for $\pi_{\nu}(x,y)$ in short ranges, with good uniformity in $\nu$. As a pleasant byproduct of these bounds, we also obtain mean estimates for the $k$-fold divisor function $\tau_k$ over short intervals which are uniform over large $k$.  
\begin{thm}\label{mainprop}
Let $\varepsilon>0$ and $a>4$ be fixed. Given sufficiently large $x \geq 1$,  let  $\mathcal{L}_{a}(x)=\log x/(\log_2 x)^a $ and suppose that $x^{17/30 + \varepsilon} \leq y \leq x$. Then we have the estimates
\begin{equation}\label{genupperlowerpinuxy}
\frac{y \delta_{\nu}(x)}{ 8^{\nu /L_{\nu}(x)} } \ll \pi_{\nu}(x,y)  \ll (1+o(1))^{\nu} y \delta_{\nu}(x)
\end{equation}
uniformly in $1 \le \nu \le \mathcal{L}_a(x)$. 

Moreover,  when $1 \leq \nu \leq (\log x)^{1/3}/(\log_2 x)^2$, we have the comparative asymptotic formula
\begin{equation}\label{shortHR}
\pi_{\nu}(x,y) \sim y \delta_{\nu}(x).
\end{equation}
Here $\delta_{\nu}(x)$ is the long-range density function defined implicitly in \eqref{HilTen}. 

\end{thm}

\begin{cor}\label{divcor}
Under the same assumptions and notation given in the first part of Theorem \ref{mainprop}, we have the estimate
\begin{equation}\label{shortdiv}
 \sum_{\substack{ n \in \halfopen{x}{x+y} \\ \Omega(n)\leq \mathcal{L}_a(x)} } \tau_k(n) \ll_{\varepsilon} y  (\log(Bk))^{11k} \ (\log x)^{2 + (1+o(1))k}
\end{equation} 
for any $k \ge 2$ and some absolute constant $B>0$. Moreover, under the additional assumption that $ 2 \leq k \leq \mathcal{L}_{\gamma}(x)$ for some $\gamma > a+1 >5$, we have the sharper bound
\begin{equation}\label{sharpshortdiv}
 \sum_{ \substack{ n \in \halfopen{x}{x+y} \\ \Omega(n)\leq \mathcal{L}_a(x) }} \tau_k(n) \ll_{\varepsilon} y   (\log k)^{11k}  \exp \Big(  \frac{\gamma +\varepsilon}{\gamma -1} k L_k(x) \Big)
\end{equation}
for any $\varepsilon > 0$.
\end{cor}

\begin{rmk}
a) In the case $\nu =1$, the stated range of validity for $y$ in the asymptotic \eqref{shortHR}  corresponds to that of Guth and Maynard's prime number theorem.\\   
b) Regarding the divisor sums in \eqref{shortdiv}, Nair and Tenenbaum \cite{NT}, and later Henriot \cite{Henriot},  obtained  short-interval estimates for a broad class of positive sub-multiplicative functions. Further, Sankaranarayanan and Srinivas \cite[Theorem 4]{SankSrini} gave an asymptotic result for short sums $\sum_{n \in [x,x+y] } \tau_k(n)$ with $y \geq x^{7/12 + \varepsilon}$,  provided that $k=\exp(o(\log_2 x))$. \\
c) Conditionally,  we believe that the arguments in this paper can easily be modified to recover \eqref{shortHR} with "full" uniformity $\nu \leq \mathcal{L}_2(x)$.  For instance, any Quasi Riemann hypothesis RH($\epsilon$) for $\zeta(s)$, which is to say a zero-free region of the form $\left\{ 1-\epsilon <Re(s) < 1 \right\}$,  should be sufficient.\\
d) Besides being a natural question in its own right, the estimation of $\| \tau_k \|_{L^{1}([x,x+y]) }$  played an important role in recent work of the author and Nishry \cite{BN}, where short divisor sums such as \eqref{shortdiv} were employed to obtain supremal bounds for random Dirichlet polynomials with Steinhaus  multiplicative coefficients.
\end{rmk}

\subsection{Notation}
$\N^{*}$ will denote the set of strictly positive integers. We will write $f \ll g$ or alternatively $f=O(g)$, if there exists an absolute constant $C$ such that $|f|\leq C |g|$. Often times we will add a subscript $f \ll_t g$ to emphasize the dependence of the implicit constant $C$ on the parameter $t$. We will use the shorthand $\log_j x$ for the $j$-fold iterated logarithm and $\tau_k$ for the $k$-fold divisor function.  A summation $\sum_{n \sim x}$ runs over integers $n \in \halfopen{x}{2x}$. We let $\omega(n)$ resp. $\Omega(n)$ denote the number of prime divisors of $n$, counted without resp. with multiplicity. $P^{-}(n)$ and $P^{+}(n)$ denote the smallest resp. largest prime divisor of $n$. Apart from section \ref{deltasection}, the symbol $\mu$ is reserved for the M\"obius function. The superscript $\flat$ will indicate a restriction to squarefree variables while the symbol $\square$ indicates perfect square integers.  Further, the set of \textsl{powerful numbers} will be denoted by $\mathscr{P}$ while the symbol $\mathcal{P}$ refers to the set of primes.  

We will make frequent use of Mertens' (second) Theorem which asserts that
\[ \sum_{p \leq x} \frac{1}{p} = \log_2 x + c+ o(1), \qquad \text{ as } x\rightarrow \infty.  \]
The parameter $\nu$ will be used to denote the number of distinct divisors of a given positive integer $n \leq x$, i.e. $\omega(n)=\nu$; we will always assume that $\nu \leq \log x$. \\

{\bf Rough subsets of $S_{v}$.} For any $1 \leq t \leq x$,  we define the set
\begin{equation}\label{lowerdeltatdef}
 \mathcal{A}_{v}(x;t)=\left\{n \in S_{v}(x): P^{-}(n) > t    \right\},   \qquad \underline{\delta}_{v} (x;t)=\frac{1}{x} \left|  \mathcal{A}_{v}(x;t)\right|.
\end{equation}

{\bf Notation for complex numbers and contour integrals.} We will use the standard notation $s=\sigma+it$ for a complex variable and, given any $c \in \R$ we will write $\int_{(c)} =\int_{c-i \infty}^{c+i \infty}$ to denote a vertical contour integral along the line $\left\{ c+it: t \in \R \right\}$. \\

{\bf A zero-free region for $\zeta(s)$.} Let us recall the celebrated Vinogradov-Korobov Theorem \cite[Section 6.19]{Tit} which ensures that the Riemann zeta function $\zeta(s)$ has no zeroes in the region
\begin{equation}\label{VinogradovKorobov}
\left\{s=\sigma+ it: \sigma > 1- \frac{C}{(\log (2+|t|))^{2/3} (\log_2 (2+|t|) )^{2/3} } \right\},
\end{equation}
for some absolute constant $C>0$.

Finally, for the convenience of the reader, we give the following short list of expressions, as they will be used on multiple occasions:
\[ \nu'_x=\frac{\log x}{\nu}, \qquad \  \ L_{\nu}(x)= \log \Big( \frac{\log x}{ \nu \log (\nu +1)} \Big).  \]
\[ \lambda^{+}=\frac{ \log x}{ \sqrt{\log_3 x} }, \qquad \tau=\exp(\lambda^{+}). \] 

\[\mathcal{L}_{a}(x)=\log x/(\log_2 x)^a , \qquad a>0,  \ x\geq 10. \] 

\[ z^{*}_{v}(c,x)=\exp\Big[ \frac{v_x'}{2 \exp(L_{v}(x)^c ) } \Big], \qquad c \in (0,1).   \]

\[\ell_{\nu}(x) = \frac{ \nu \ \log(L_{\nu}(x))^2 }{L_{\nu}(x) }. \]

\subsection{Acknowledgements} This work was supported by the  Academy of Finland grant CoE FiRST. Part of the research was carried out during the 2024 thematic semester on analytic number theory at the Mittag-Leffler institute; I gratefully acknowledge their hospitality and support. I am indebted to Neea Paloj\"arvi and Ofir Gorodetsky for many fruitful discussions, helpful comments and references. I am also very grateful to the anonymous referee for their careful reading of the article,  their helpful suggestions and comments.

\section{Preliminaries}
\subsection{The density function \boldmath\texorpdfstring{$\delta_{\nu}$}{delta{nu}}.} \label{deltasection}As briefly touched upon in the introduction, Hildebrand and Tenenbaum \cite{HT} gave an asymptotic formula for the cardinality of $S_{\nu}(x)$ which is valid in the range $\nu \leq \log x/ (\log_2x )^2=\mathcal{L}_2(x)$. The purpose of this first section is to give a simple expression for the (approximate) order of $\delta_{\nu}$, allowing for a direct comparison with Landau's density function given in \eqref{Lan}. 

To this end, rather than working with the general expression \eqref{HilTen}, it will be convenient to separate those values of $\nu $ below and beyond the threshold $(\log_2 x)^2$.  Accordingly, we may use \cite[Main theorem, eq. (2.7)]{HT} for smaller values of $\nu$ and \cite[Corollary 2]{HT} for large values. We have that \\
\begin{equation}\label{deltaformula}
\delta_{\nu}(x) \sim
\left\{
	\begin{array}{lll}
		\frac{1}{\nu ! (\log x)} \exp(\nu(\log M +\frac{1}{M} +O(R) ))  & \mbox{if }  (\log_2 x)^2 \leq \nu \leq \mathcal{L}_2(x) \\
		& \\
		\frac{(\log_2 x)^{\nu }}{ \nu ! \log x} \rho H(\rho) \exp(\nu ( \gamma- \log(1+ \gamma)  ) ) & \mbox{if } \nu< (\log_2 x)^2.
	\end{array}
\right.
\end{equation}
Here $H$ is the Euler product
\begin{equation*}
H(s)=\frac{1}{\Gamma(s+1)} \prod_{p} \Big( 1+ \frac{s}{p-1}  \Big) \Big(1-\frac{1}{p} \Big)^s
\end{equation*}
and $\rho$ is the first coordinate of the unique minimising pair \footnote{The existence and uniqueness of $(\rho, \alpha)$ was shown in \cite[Lemma 2]{HT}.}
\begin{equation}\label{minimisers}
 (\rho, \alpha)=(\rho_{v}(x), \alpha_v(x))= \mathop{{\arg \inf}}_{(r,a) \in (0,\infty) \times (1,\infty)} G(r,a)x^{a} r^{-\nu}, 
 \end{equation}
where $G$ is the generating series defined in \eqref{Gdef}. Moreover, the quantities $L=L_{\nu}(x),M,R$ are given by
\begin{equation}\label{HTparameters1}
L= \log \Big( \frac{\log x}{\nu \log(\nu +1)} \Big), \qquad M=\log \Big(\frac{C w \log w}{L}   \Big), \qquad R=\frac{1}{L \log(\mu +2)} +\frac{1}{L^2} 
\end{equation} 
with
\begin{equation}\label{HTparameters2}
u= \frac{\nu}{\log_2 x}, \qquad  \mu=\frac{\nu}{L_{\nu}(x) }, \qquad w=\frac{\log x}{\mu \log(\mu +2)}, \qquad \gamma=(\rho - u)/u 
\end{equation}
and $C>0$ being an absolute constant. In order to work with the implicitly defined minimising pair $(\rho_{v}(x), \alpha_v(x))$, it will be helpful to recall their approximate orders, which were determined in \cite[Lemma 2]{HT}: for $\nu \leq \mathcal{L}_2(x)$ we have that
\begin{equation}
   \rho= \rho_{\nu}(x) = \frac{\nu}{L} + O\Big( \frac{\nu \log L}{L^2} \Big), \qquad \alpha_v(x) =1+ \frac{\rho}{ \log x} \Big(1+O \Big( \frac{\rho \log(\rho +10)}{\log x} \Big) \Big).
\end{equation}

It should be noted that the first line of \eqref{deltaformula} does not yield an asymptotic formula for $\delta_{\nu}(x)$. Instead, we will use the above expressions to derive crude upper and lower bounds for the density function.  
\begin{lem}\label{4divisor}
With notation as above, one has upper and lower bounds of the form 
\begin{equation}\label{deltacrude}
 \delta_{\nu}(x)  = \frac{( e_{\nu}(x) L_{\nu}(x) )^{\nu  -1}}{( \nu -1) ! \log x}. 
\end{equation}
The functions $e_{\nu}(x)$ are uniformly close to $1$ in the sense that $e_{\nu}(x)=1+O(\eta(x))$ for some $\eta(x)=o(1) $ and all $\nu \leq \mathcal{L}_2(x)$.  

\begin{proof}
First, when $ (\log_2 x)^2 \leq \nu \leq \mathcal{L}_2(x)$, we require bounds for the quantities $L, M,R$ given in \eqref{HTparameters1} and \eqref{HTparameters2}. In this range of $\nu$, we see that $L \in [\log_3 x, \log_2 x]$ for all large $x \geq 1$ and hence $\log \mu \in [\frac{1}{2} \log \nu, \log \nu ]$. We find that
\[ \log w= \log \Big( \frac{L \log x}{\nu \log(\mu + 2) } \Big) \sim \log L + \log \Big( \frac{ \log x}{\nu \log(\nu + 1) } \Big) \sim L  \]
as $x \rightarrow \infty$, from which we get that $M \sim \log w \sim L$ and $R=o(1)$. It follows that the last two terms appearing inside the exponential in the first line of \eqref{deltaformula}  are bounded by $|M^{-1} + O(R)| =o(1)$ and hence we retrieve \eqref{deltacrude} upon noting that $L_{\nu}(x)/ \nu =(1+o(1))^{\nu}  $.\\
\indent In the short range $\nu \leq (\log_2x)^2$, we make use of the approximation $\rho = \nu /L + O(\nu \log L /L^2)$ and observe that $L= \log_2 x + O(\log_3 x)$. From \eqref{HTparameters2} we gather that  $1+ \gamma =\rho /u=1+ O(\log_3 x/ \log_2 x)$, whence $ \gamma- \log(1+ \gamma  )= O(\log_3 x/ \log_2 x) $. Since $\rho \log_2 x / \nu \in [1/2,2]$, all that remains is to estimate $H(\rho)$. First, under the assumption that $\rho >100$ we find
\begin{align*}
 &C_{\rho}:= \log (\Gamma(\rho +1)  H(\rho)) = \sum_{p}  \big[ \log(1+ \frac{\rho}{p-1}) + \rho \log (1- \frac{1}{p}) \big]  \\
 & =\sum_{p \leq 2 \rho}  \big[ \log(1+ \frac{\rho}{p-1}) - \frac{\rho}{p}) \big]  +O \Big(\rho+ \sum_{p > 2\rho }  \frac{\rho^2}{p(p-1)}  \Big) \\
 &=  - \rho \log_2 \rho + O(\rho).
 \end{align*}
An application of Stirling's formula yields
\[H(\rho)= \frac{ \exp(C_{\rho}) }{  \Gamma(\rho +1)}  =\frac{ \exp( O(\rho)) }{ \sqrt{\rho}( \rho  \log \rho )^{\rho}}=\exp( o(\nu))    , \]
which may be combined with our previous estimate $ \gamma- \log(1+ \gamma  )= o(1)$ and inserted into the second line of \eqref{deltaformula} to get \eqref{deltacrude}. 

Finally, when $\rho \leq 100$, our treatment of $C_{\rho}$ reveals that $H(\rho) \asymp 1= (1+o(1))^{\nu} $ in the range $\sqrt{\log_2 x} \leq \nu \leq (\log_2 x)^2$, whereas $H(\rho) =1+o(1)$ when $\nu <\sqrt{ \log_2 x}$. In either case \eqref{deltacrude} remains valid.
\end{proof}
\end{lem} 

To conclude the discussion of $\delta_{\nu}$,  we record some useful comparative statements. First we have that \cite[Corollary 4]{HT}    
\begin{equation}\label{deltahomothety}
\delta_{v}(mx)= \delta_{v}(x) \Big( \frac{\log(mx)}{\log x}\Big)^{\frac{v}{L} -1} \exp\Big( O\big (\frac{1}{L} + \frac{v \log L \log m }{L^2 \log x} \big) \Big), 
\end{equation}
with $L=L_v(x)$.  This approximate identity/estimate holds for all large $x$,  and all $v \leq \mathcal{L}_2(x)$ and $1\leq m \leq x$.  The second result, given in \cite[Corollary 3]{HT} asserts that
\begin{equation}\label{densitycomparisonv}
\frac{\delta_{v+1}(x) }{\delta_{v}(x) } \sim \frac{L_v(x) }{v}  , \qquad  v \leq \mathcal{L}_2(x). 
\end{equation}

\subsection{Counting integers in \texorpdfstring{$S_{\nu}$}{S{nu}} with multiplicative constraints}

In this section we study various subsets of $S_{v}$ that are made up of integers whose prime divisors lie in prescribed ranges. Given an interval $I \subset [1,x]$ and $v \ge 1$, we first set
\begin{equation}
\omega_{I}(n)=\sum_{p \in I, p|n} 1, \qquad S_{I,v}(x)=\left\{ n \in S_{v}(x): \omega_{I}(n)=v \right\}.
\end{equation}
Next we describe a type of exceptional set that will be used extensively in the proof of \eqref{shortintervalminorantupper}.  For $v \leq \mathcal{L}_2(x)$ and any $c \in (0,1)$, we first recall the notation $v'_x=(\log x)/ v$ and set
\[ \log z_v^{*}(c,x)=  \frac{ v'_x }{ 2\exp( L_v(x)^{c} )} , \qquad  I^{*}_{v}(c,x)=\big[z^{*}_{v}(c,x),  z^{*}_{v} \big(\frac{c}{2},x \big) \big].  \]
Sometimes we will suppress the dependence on $x$ and simply write $z_v^{*}(c)$ or even $z^{*}_v$ to simplify notation.  Throughout the remainder of this article, we will almost always assume that $x$ is large,  $v \leq \mathcal{L}_2(x)$ and that
\begin{equation}\label{widecrange}
\frac{1}{\sqrt{\log(L_v(x) )} } \leq c < 9/10,
\end{equation}
in which case we have the straightforward bounds
\begin{equation}\label{solnapprox}
    v+1 \leq z_v^{*}(c) \leq x^{o(1)}.
\end{equation}
Given real parameters $1 \leq t <x \leq x'$ , along with a partition $v=v_0+v_1$, we may now proceed to define the exceptional set
\[ S_{v,v_1, x'}^{*}(x; c;t)=\left\{n \in S_{v}(x): \omega_{I^{*}_{v_1}(c,x') }(n)=v_1,  \ P^{-}(n) > t    \right\},   \]
which is made up of $t$-rough integers $n \in S_v(x)$ having exactly $v_1$ distinct prime divisors in the interval $I^{*}_{v_1}(c,x')$. Its natural density will be denoted by
\[  \delta_{v,v_1,x'}^{*}(x; c;t)= \frac{1}{x}|S_{v,v_1,x'}^{*}(x;c;t)|,\]
and in the special case where $x=x'$ and $t=1$, we write  $S_{v,v_1}^{*}(x; c)$ and $ \delta_{v,v_1}^{*}(x; c)$. Further, for any $u \leq v$, we set
\begin{equation}\label{Deltastardef}
 \Delta^{*}_{v,u,x'}(x; c;t) = \sum_{u \leq w \leq v}   \delta_{v,w,x'}^{*}(  x;c;t).
\end{equation}

\begin{lem}
For any interval $I= [a,b] \subset [2,x]$ and any $v \geq 1$, we have the estimate 
\begin{equation}\label{Fnuprescribed}
\sum_{n \in S_{I, v }(x) } \frac{1}{n} \ll  \frac{(\log_2 b -\log_2 a +\varepsilon(x)+\varepsilon(a)  )^{v}}{v !}. 
\end{equation}
Under the additional assumption that $v \leq \mathcal{L}_1(x)$, we have the bound
\begin{equation}\label{Fnuprescribednatural}
|S_{I, v }(x)|   \ll  \frac{x v}{\log x} \frac{(\log_2 b -\log_2 a +\varepsilon(x)+\varepsilon(a)  )^{v-1}}{(v-1) !}. 
\end{equation}
Here the bounded function $\varepsilon(y)$ tends to zero as $y \rightarrow \infty$. Moreover, given a non-trivial partition $v=v_0+v_1 $, any $c>0 $ as in \eqref{widecrange} and parameter $1 \leq t \leq (z_v^{*}(c,x) )^{1/2} $, we have the lower bound 
\begin{equation} \label{lowerrestrictedstar}
\delta_{v,v_1,x'}^{*}(x;c;t) \geq \frac{1}{2^v v !  \log x} \Big( \log \Big(\frac{\log z_v^{*}(c,x) }{\log 2t}  \Big)  \Big)^{v_0}  L_{v_1}(x) ^{cv_1} 
\end{equation}
for all sufficiently large $x\geq1$, all $v \leq \mathcal{L}_2(x)$ and any $x' \in [x,x^{3/2}]$.
\begin{proof}
The estimate \eqref{Fnuprescribed} is a straightforward consequence of the combinatorial observation
\[ \sum_{n \in S_{I, v }(x) } \frac{1}{n}
\leq \frac{1}{v !} \Big( \sum_{p \in I}  \sum_{k \geq 1} \frac{1}{p^k} \Big)^{v}, \]
followed by an application of Mertens' Theorem. The bound \eqref{Fnuprescribednatural} follows immediately from \eqref{Fnuprescribed} since, outside an exceptional subset $\mathscr{E} \subset S_{I, v }(x)$ of size $O(x^{3/4})$, every $n \in S_{I, v }(x)$ satisfies $P^{+}(n)\geq x^{1/2 v}$ and thus
\[ |S_{I, v }(x)| \leq  |\mathscr{E}|+ \sum_{m \in S_{I, v-1 }(x^{1-1/2v}) } \ \sum_{ \substack{ p^e \leq x/m \\ e\geq 1}} 1 \ll x^{3/4} + \frac{x v}{\log x}  \sum_{m \in S_{I, v-1 }(x) } \frac{1}{m}.  \] 
 Next we address the lower bound \eqref{lowerrestrictedstar}. Observe that by choosing $v_1$ distinct primes in the interval $I_{v_1}^{*}(c,x') $, $v_{0}-1$ primes in the range $[2t, z^{*}_v(c,x') ]$ and one large prime, we get the lower bound
\begin{align*}
|S_{v,v_1,x'}^{*}(x; c;t)|   &\geq  \frac{1}{ v_0 ! v_1!}  \sum_{ \substack{ (p_1,...,p_{v_1}) \in (I_{v_1}^{*}(c,x') )^{v_1} \\  n=p_1\cdots p_{v_1} : \mu^2(n)=1 }} \ \ \sum_{ \substack{ (p'_1,...,p'_{v_{0} -1 })  \\  m=p'_1\cdots p'_{v_{0}-1} : \mu^2(m)=1 \\ 2t \leq p_j' \leq z_{v}^{*}(c,x') } } 
\sum_{\sqrt{x} < p \leq \frac{x}{mn} } 1 
\\
\gg  \frac{x}{(\log x)v_0 ! v_1!}    & \sum_{ \substack{ (p_1,...,p_{v_1}) \in (I_{v_1}^{*}(c,x') )^{v_1} \\  n=p_1\cdots p_{v_1} : \mu^2(n)=1  }}     \frac{1}{p_1 \cdots p_{v_1}}       
\sum_{ \substack{ (p'_1,...,p'_{v_{0} -1 })  \\  m=p'_1\cdots p'_{v_{0}-1} : \mu^2(m)=1 \\ 2t \leq p_j' \leq z_{v}^{*}(c,x') } }      \frac{1}{p'_1 \cdots p'_{v_{0}-1 }}
\end{align*}
from which we may retrieve  \eqref{lowerrestrictedstar} after repeated application of Mertens' Theorem.
\end{proof}
\end{lem}

Our next goal is to show that, under suitable conditions, the exceptional set $S_{v,v_1}^{*}(x; c)$ is indeed small, in which case we wish to produce a sufficiently sharp upper bound for its cardinality $\delta_{v,v_1}^{*}(x;c)$. In Lemma \ref{deltastarlemma} just below, we will give an estimate for the cardinality $\delta_{v,v_1}^{*}(x;c)$ that is sensitive to both "overcrowding" and sparsity of prime divisors in the interval $ I_{v_1}^{*}(c,x)$. Before moving on to the treatment of $S_{v,v_1}^{*}(x; c)$ we define, for any $c \in (0,1)$, the modified Euler product
\[ G^{*}_{v,c,x}(z,s) =\prod_{p \in I_v^{*}(c,x)  }  \Big( 1+  \frac{z}{p^s -1}\Big). \]

We recall that $(\rho_{v}(x), \alpha_v(x))$ is the minimising pair given in \eqref{minimisers}.
\begin{lem}\label{deltastarlemma}
Let $x\geq 1$ be sufficiently large and assume that $1 \leq v \leq \mathcal{L}_2(x)$.  Given a coefficient $c$ as in \eqref{widecrange} and a partition $v=v_0+v_1$,  satisfying
\begin{equation}\label{vonecondition}
 \frac{v}{(L_v(x))^{10}} \leq v_1 \leq v,  
\end{equation}
we have that
\begin{equation}\label{deltastarmainbound} 
\delta_{v,v_1}^{*}(x;c)\ll   \delta_{v}(x)   v_1^2 (\log_2 x)^3  \frac{(3 \gamma)^{v_1} }{\exp(\gamma v_1 /2)},
\end{equation} 
where $\gamma=\gamma_{v,v_1}(c)=\rho_v(x) L_{v_1}(x)^c / v_1$. 
\begin{proof}
Let $x^{1/10} \leq y \leq x$ and write $ I_{v_1}^{*}(c)= I_{v_1}^{*}(c,x)$ to ease notation. Since any $m \in S_{v_0}(y)\cap (y/2,y] $ (outside a small exceptional set) has a prime divisor $p \geq x^{1/20 v_0}$ we first find, as in the proof of \eqref{Fnuprescribednatural}, that
\[   \sum_{ \substack{m \in S_{v_0}(y)\cap (y/2,y]     \\ p|m \implies p \notin  I_{v_1}^{*}(c) }   } 1 \ll \frac{v_0}{\log x} \ y \sum_{ \substack{n \in S_{v_0 -1}(y)    \\ p|n \implies p \notin  I_{v_1}^{*}(c) }   } \frac{1}{n}.  \]
An application of Rankin's trick, together with the estimate $v/\rho_v \ll \log_2 x$,  yields

\begin{align}\label{Rankinapproach}
   \sum_{ \substack{m \in S_{v_0}(y)\cap (y/2,y]     \\ p|m \implies p \notin  I_{v_1}^{*}(c) }   } 1&  \leq \frac{v_0}{\log x} \sum_{\substack{ m \geq 1 \\ p|m \implies p \notin  I_{v_1}^{*}(c) } }  \rho_{v}^{\omega(m)-v_0 -1}
   \big(\frac{y}{m} \big)^{\alpha_{v} }       \notag   \\
& \ll \frac{ \log_2 x}{\log x} y^{\alpha_{v} } \rho_{v}^{-v_0}  \frac{G(\rho_{v},\alpha_{v})}{  G_{v_1,c,x}^{*}(\rho_{v},\alpha_{v})}  \leq y  \frac{ \log_2 x}{\log x} \ \rho_{v}^{-v_0}  \frac{x^{\alpha_{v} -1}  G(\rho_{v},\alpha_{v})}{  G_{v_1,c,x}^{*}(\rho_{v},\alpha_{v})}   \\
& \sim y \delta_{v}(x) \big( v \phi(v) \phi(\rho_{v}) \log_2 x \big) 
\frac{  \rho_v^{v_1}}{  G_{v_1,c,x}^{*}(\rho_{v},\alpha_{v})} ,  \notag
\end{align}
where $\phi$ was defined in \eqref{psiLdef}. By Stirling's approximation we have that $\phi(y) \ll y^{-1/2}$ for $y\geq 1$ and hence, recalling the asymptotic $\rho_{v}(x) \sim v/ L_{v}(x) $, we find that $v \phi(v) \phi(\rho_{v}) \log_2 x \ll (\log_2 x)^2$. Turning to the Euler products, we combine the approximation $\alpha_v \sim 1+ \rho_v/ \log x$ with \eqref{vonecondition} to get that $p^{-\alpha_v}\sim 1/p$ for any $p \in I_v^{*}(c)$.  From the lower bound in \eqref{solnapprox} it follows that
\[ \log  G^{*}_{v_1,c,x}(\rho_v,\alpha_{v} ) = \rho_v \sum_{p  \in I_{v_1}^{*}(c)  }  p^{-\alpha_v} + O(\rho_v^2 / v) =  \rho_v(x) L_{v_1}(x)^c(1+o(1)) .  \]

Observe that any $n \in S_{I_{v_1}^{*}(c), v_1} $ is of size at least $z^{*}_{v_1}(c,x)^{v_1}$ and at most $\sqrt{x}$ (since each prime divisor $p|n$ is smaller than $x^{1/2v_1}$, say).  Writing $u=x/z^{*}_{v_1}(c,x)^{v_1}$, it follows from \eqref{Fnuprescribednatural} and \eqref{Rankinapproach} (applied to each dyadic interval between $\sqrt{x}$ and $x$) that
\begin{align*}
   \delta_{v,v_1}^{*}(x;c) & \leq x^{-1} \sum_{ \substack{m \in S_{v_0}(x) \\  \sqrt{x} \leq m \leq  u   \\ p|m \Rightarrow p \notin  I_{v_1}^{*}(c) }   } S_{I_{v_1}^{*}(c),v_1 } \big(\frac{x}{m} \big)  
 \ll   \sum_{ \substack{m \in S_{v_0}(x),m\geq \sqrt{x}     \\ p|m \Rightarrow p \notin  I_{v_1}^{*}(c) }   }  \frac{v_1 \big(  L_{v_1}(x)^c [1+o(1)] \big)^{v_1 -1}  }{ (v_1 -1) ! \ m \ \log x}  \notag \\
 & \ll \delta_{v}(x)  v_1^2 (\log_2 x)^3  \ \frac{\rho_v^{v_1}  \big( L_{v_1}(x)^c [1+o(1)] \big)^{v_1}  }{  G_{v_1,c,x}^{*}(\rho_{v},\alpha_{v}) \ v_1 !}  \\
 &\ll \delta_{v}(x)   v_1^2 (\log_2 x)^3   \frac{(3 \gamma)^{v_1} }{\exp(\gamma v_1 /2)},
 \end{align*}
as desired. 
\end{proof}
\end{lem}

In the following lemma we are given a partition $v =v_0 +v_1$ together with subsets $\mathcal{C}_0 \subset S_{v_0}$ and  $\mathcal{C}_1 \subset S_{v_1}$.  Under suitable multiplicative constraints $\mathcal{C}_0 $ and $\mathcal{C}_1 $ may be reassembled or "glued together"  to form a subset of $S_{v}$.  

\begin{lem}\label{gluglu}  Let $\eta \in (0,1)$ be fixed, suppose that $x\geq 10$ is sufficiently large and let $v \leq \mathcal{L}_2(x)$.  Further, let $v=v_0 +v_1$ be a partition satisfying \eqref{vonecondition} and $v_1 \geq 1/\eta$.  Given any coefficient $c>0$ satisfying \eqref{widecrange} and any $1 \leq t \leq z_v^{*}(c,x) $,  we have the estimate
\begin{equation}\label{Gluing}
 \sum_{\substack{n \in S_{v_0 }(x) \\ P^{+}(n) \leq t }}  \frac{1}{n}  \Delta^{*}_{v_1, \eta v_1,x} \big( \frac{x}{n};c;t\big) \ll \delta_{v}(x)  v_{1}^2 (\log_2 x)^3   \frac{ (3 \gamma/ \eta)^{v_1}  + (3 \gamma/ \eta)^{\eta  v_1}  }{\exp(\gamma v_1 /3)},
  \end{equation}
where $\gamma=\gamma_{v,v_1}(c)$ is the discriminant defined in Lemma \ref{deltastarlemma}. 
\begin{proof}

Observe that every $n \in S_v$ admits at most one representation $n=n_0n_1$ satisfying both $P^{+}(n_0) \leq t$ and $P^{-}(n_1) > t$  . As a consequence, we have that
\begin{align*} 
 \sum_{\substack{n \in S_{v_0 }(x) \\ P^{+}(n) \leq t }}  \frac{x}{n}  & \Delta^{*}_{v_1, \eta v_{1},x} \big( \frac{x}{n};c;t\big)  = \sum_{w \in [\eta v_{1},v_1 ]}    \sum_{\substack{n \in S_{v_0 }(x) \\ P^{+}(n) \leq t }}   |S_{v_1,w,x}^{*} \big ( \frac{x}{n} ;c;t \big)|  \\
& \leq   \sum_{w \in [\eta v_{1},v_1 ]}  |S_{v,w}^{*}(x;c)| =  x \sum_{w \in [\eta v_{1},v_1 ]}  \delta_{v,w}^{*}(x; c) 
\end{align*}

and hence, taking into account the relationship  $  L_{w}(x) =(1+o(1)) L_{v_1}(x) $ (for $w \in[\eta v_{1},v_1]$), the result follows from Lemma \ref{deltastarlemma}.
\end{proof}

\end{lem}

\subsection{ \boldmath\texorpdfstring{$S_{\nu}$}{S{nu}} in short intervals}

Following the discussion in section \ref{priorsection}, the next step is to move from small values of $\nu$ to the case of slowly growing $\nu$.

\subsubsection{The range $\nu=o( \log_2 x) $} To give a complete picture of the range $\nu \leq \log_2 x /(\log_3 x)^{3/4}$,  we will make use of a result of K\'atai.  Inserting the zero density estimate of Guth-Maynard (see \eqref{Hux} below)  into \cite[Main Theorem]{Kat} and applying \eqref{deltaformula}, one immediately finds that

\begin{equation}\label{Avformula}
\pi_{\nu}(x,y) \sim y \delta_{\nu}(x) \sim \frac{y}{\log x} \frac{(\log_2 x)^{\nu -1}}{(\nu -1)!}=: y s_{\nu}(x), \qquad \nu \leq \frac{\log_2 x}{(\log_3 x)^{3/4}},
\end{equation} 
uniformly in $x^{17/30 + \varepsilon} \leq y \leq x$.  In other words,  the short-range density of $S_{\nu}(x,y)$ coincides with Landau's density $s_{\nu}(x)$ appearing in \eqref{Lan}. As a consequence, we obtain the anatomical statement $\pi_{\nu}(x,y) \sim \mathcal{M}_{\nu}'(x,y)$.\\

\noindent {\bf Proof of \eqref{simpleanatomy}. } Fix $\varepsilon>0$, and let $x^{17/30 + \varepsilon} \leq y \leq x$. It is straighfoward to verify that $s_{\nu}(\tau) \sim s_{\nu}(x)$ and $s_{\nu -1} (\tau) =o(s_{\nu}(\tau))$ for any $\nu \leq \log_2 x / (\log_3 x)^{3/4}$, and thus it follows from the Guth-Maynard prime number theorem \cite[Corollary 1.3]{GM} that
\begin{align*}
 \mathcal{M}_{\nu}' (x,y)& =  \sum_{m \in S_{\nu -1}(\tau) }  \sum_{p \in \halfopen{x/m}{(x+y)/m}}   \sim \frac{y}{\log x} \sum_{m \in S_{\nu -1}(\tau) }  \frac{1}{m} \\
 &= \frac{y}{\log x} \int_{1}^{ \tau} \frac{d \pi_{\nu -1}(u) }{u}   = \frac{y}{\log x}  \left[ \delta_{\nu -1} (\tau) + \int_{1}^{\tau} \frac{\delta_{\nu-1}(u)  }{u} \ du\right] \\
& \sim \frac{y}{\log x}  \left[ s_{\nu -1} (\tau) + \int_{1}^{\tau} \frac{s_{\nu-1}(u)  }{u} \ du\right] \sim y s_{\nu}(\tau) \sim \pi_{\nu}(x,y).
 \end{align*}

\begin{rmk}
It is worth noting that the anatomical property $\pi_{\nu}(x,y) \sim \mathcal{M}_{\nu}'(x,y)$ is no longer valid when $\nu \gg \log_2 x$.  Indeed, we have that
\begin{equation*}
\sum_{x^{1/4} \leq m \leq \sqrt{x} } \ \sum_{p \leq x/m} \gg x
\end{equation*}
and thus a positive proportion of the integers $n \leq x$ are not expressible as a product of the form $mp$ with $m \leq \tau$. 
\end{rmk}

\subsubsection{A uniform lower bound for $\pi_{\nu}(x,y)$} Our final task is to produce a lower bound for $\pi_{\nu}(x,y)$ in the long range $\nu \leq \mathcal{L}_4(x)$. We will make use of Wu's result \cite[Theorem 1  ]{Wu} which gives a lower bound for the number of $2$-almost primes in short intervals:
\begin{equation}\label{twoalmostlower}
 \sum_{\substack{ m \in \mathcal{P}_2 \\ m \in \halfopen{x}{x + x^{\varpi} } }} 1 \geq c \frac{x^{\varpi}}{\log x}, \qquad \ \varpi=101/232 \approx 0.4354
\end{equation}
for some constant $c>0$ and all $x \geq 2$.  Here, $\mathcal{P}_2 \subset S_1 \cup S_2$ denotes the set of integers with at most two prime factors, counting multiplicity. \\

{\bf Proof of the lower bound \eqref{genupperlowerpinuxy}, assuming \eqref{shortHR}.}
We set $H=x^{1/7.9 }$ and assume that  $\nu \geq (\log_2 x)^2$.   From \eqref{deltahomothety} and \eqref{densitycomparisonv} we first gather that
\[ \min(\delta_{\nu -2}(H),  \delta_{\nu -1}(H)) \gg  (7.9)^{-\nu / L_{\nu}(x) } \delta_{\nu}(x). \] 
An application of \eqref{twoalmostlower} yields
\begin{align*} 
 \pi_{\nu}(x,y) & \gg \frac{1}{\nu^2} \sum_{j\in \left\{ 1, 2 \right\}  } \big(  \sum_{m \in S_{\nu -j}(H) } \ \sum_{\substack{ n \in \mathcal{P}_j \\ n \in \halfopen{x/m}{ (x+y)/m  } }} 1  \big) \\
  &\gg \frac{y}{\nu^2 \log x}  \ \min_{j\in \left\{ 1, 2 \right\}  } \big( \sum_{m \in S_{\nu -j}(H) } \frac{1}{m} \big) \gg \frac{y}{\nu^2 \log x} \min_{j\in \left\{ 1, 2 \right\}  }  \delta_{\nu -j}(H)
  \end{align*}  
and thus we recover \eqref{genupperlowerpinuxy}.

\section{The bilinear and unweighted sums \texorpdfstring{$\mathcal{F}_{\underline{\lambda}}$ and $\mathcal{W}$.}{F{lambda} W{lambda} and V{lambda} } }

\subsection{The set-up and outline of the strategy}\label{outlinesection}
Having completed the preliminary work of the previous section, we are ready to embark upon the main part of our journey.  On the one hand we are concerned with the derivation of the asymptotic formula \eqref{shortHR} asserted in Theorem \ref{mainprop}. On the other hand,  we wish to understand the anatomical claim made in Theorem \ref{squarefreemainprop}.  Regarding the latter assertion, our goal is to show that the vast majority of integers in $S_{\nu}(x,y)$ cannot have a large concentration of prime divisors of order approximately $x^{1 / \nu}$.  In other words,  the \textsl{typical integer} in $S_{\nu}(x,y)$ is overwhelming made up of "small" primes.  In order to realise this plan (which will be carried out in Section \ref{upperboundproofsection} below), we require a technical device that is able to detect and bound the size of exceptional subsets of $S_{\nu}(x,y)$ with sufficient accuracy. For this purpose we will employ the bilinear sums $\mathcal{F}_{\underline{\lambda}}$, which are constructed as follows.

Let $F_{v}$ denote the $v$-fold convolved von Mangoldt function, as defined inductively via the relations
\begin{equation}\label{Fvdef}
 F_{1}= \Lambda, \qquad   F_{v}=  F_{v -1} * \Lambda , \qquad v \geq 2. 
\end{equation} 
In particular, we see that
\begin{equation}\label{Fvsupport}
\supp(F_{v}) \subset \left( \bigcup_{w \leq v} S_{w} \right) \cap \left\{n \in \N^{*}: \Omega(n) \geq v \right\}
\end{equation}
and, moreover, when restricted to squarefree integers, the support of $F_v$ is exactly $S_v$. We may then introduce bilinear decompositions of the shape
\begin{equation}\label{curlyFdef}
 \mathcal{F}_{\underline{\lambda}}(x,y)= \sum_{ m  \sim U } \  \  \sum_{\substack{  l \in (\frac{x}{m}, \frac{x+y}{m} ] \\ P^{-}(l) >z }} F_{v}(l)
 \end{equation} 
with parameter set $\underline{\lambda}=(v,z,U)$.  

It should be noted that the choice of von Mangoldt weights is natural since the generating Dirichlet series of $F_v$ is $(-\zeta'/\zeta(s))^v$ which is a manageable function inside the critical strip, owing to a variety of well-known multiplicative and analytic techniques.  In parallel with the study of $ \mathcal{F}_{\underline{\lambda}}(x,y)$, we will examine the unweighted sum $\mathcal{W}_{\underline{\lambda}} (x,y)$,  obtained by replacing each von Mangoldt component appearing in $F_{v}$ with $\overline{\theta}=\Lambda/ \log $ -- the latter weight serves as an \textsl{almost indicator function of the primes}.  Importantly,  the unweighted sum $\mathcal{W}_{\underline{\lambda}} $ may be expressed as a suitable combination of sums which are essentially of the form $ \mathcal{F}_{\underline{\lambda}}(x,y)$ (without the additional averaging over $m \sim U$).  This fact will ultimately yield the comparative asymptotic $\mathcal{W}_{\underline{\lambda}} (x,y) \sim (y/x) \mathcal{W}_{\underline{\lambda}} (x)$.  The final step, which  will be carried out in section \ref{finalsection}, is to deduce Theorem \ref{mainprop} from this comparative result. 

Let us move on to the construction of $\mathcal{W}_{\underline{\lambda}}$.  Write $\overline{\theta}(n)=\Lambda(n)/ \log n$ and define the $k$-fold convolution
\begin{equation*} 
P_k(n)= \Big( \overline{\theta}* \cdots* \overline{\theta} \Big)(n). 
\end{equation*}
Observe that the squarefree restriction of $P_k$ picks out exactly those integers in $S_k^{\flat}$.  We may now define the long and short unweighted sums 
\begin{align}
\begin{split}
\mathcal{W}^{v,z}(x)&= \sum_{\substack{n \leq x \\ P^{-}(n) >z }} P_v(n),  \\
  \mathcal{W}^{v,z} (x,y)& = \mathcal{W}^{v,z} (x+y) - \mathcal{W}^{v,z} (x).
\end{split}
\end{align}

In the following lemma, we give a reformulation of the bilinear sum $\mathcal{F}_{\underline{\lambda}}(x,y)$ in moderately short ranges of $y$.

\begin{lem}\label{longylemma}
Let $x \geq 1$ be sufficiently large and let $v \leq \mathcal{L}_2(x)$. Suppose, moreover, that $1 \leq U \leq x/ (\log x)^{4v} $, $ U \leq y \leq x$ and that $yU\geq x(\log x)^{5v}$. Then 
\begin{equation}\label{curlyFmeanmoderateshort}
\mathcal{F}_{\underline{\lambda}}(x,y) = y\sum_{\substack{ l \in  \halfopen{ \frac{x+y}{2U} }{  \frac{x}{U} } \\ P^{-}(l) >z }} \frac{ F_{v}(l) }{l}   + O(y/ (\log x)^{2v}). 
\end{equation}

\begin{proof}
We may assume without loss of generality that $y \leq x/ (\log x)^{4 v}=: X $ (or else separate the interval $(x,x+y]$ into blocks of length $X$).  Let us interchange the order of summation in \eqref{curlyFdef} and separate the variable $l$ into the regions $I_1= \halfopen{\frac{x}{2U}}{\frac{x+y}{2U} } $, $I_2= \halfopen{ \frac{x+y}{2U} }{  \frac{x}{U} }$ and $I_3=\halfopen{\frac{x}{U} }{ \frac{x+y}{U} }$. Correspondingly, we set $A^{(1)}_l=\halfopen{x/l}{2U} $ and $A^{(3)}_l=\halfopen{U}{\frac{x+y}{l} } $ and observe that $|A^{(j)}_{l}|\leq 2Uy/x$ for any $l \in I_j$, with $j=1,3$.  An application of the straightforward pointwise bound $F_{v}(l) \leq v^v (\log x)^{2v} $ gives the desired expression

\begin{align*}
\mathcal{F}_{\underline{\lambda}}(x,y)  & = \sum_{j \in \left\{1,3 \right\} }  \sum_{\substack{  l \in I_j \\ P^{-}(l) >z }} F_{v}(l)  \sum_{m \in A^{ (j) }_{l} }   1 
+ \sum_{\substack{  l \in I_2 \\ P^{-}(l) >z }} F_{v}(l)  \sum_{m \in [\frac{x}{l}, \frac{x+y}{l} ] }  1 
\notag \\
& = O\Big( \frac{y^2}{x} v^v (\log x)^{2v} + \frac{x}{U} v^v (\log x)^{2v} \Big) + y\sum_{\substack{  l \in \halfopen{ \frac{x+y}{2U} }{  \frac{x}{U} } \\ P^{-}(l) >z }} \frac{ F_{v}(l) }{l}.  
\end{align*}

\end{proof}
\end{lem}

\subsection{Recasting the unweighted sums}
\subsubsection{Step one: some initial manipulations} As outlined in the previous section, our goal is to express the sum $\mathcal{W}^{v,z}(x,y)$ as a weighted average of sums of the form  $ \mathcal{F}_{\underline{\lambda}}(x,y)$.  To this end,  we expand the convolution sum factor-per-factor, separating each summation variable according to the dyadic ranges
\begin{equation*}
\mathscr{D}_z(x)=\left\{ 2^n z: n \in \N^{*} \cup \left\{ 0 \right\}  \right\} \cap [1,2x],
\end{equation*}
and then sum by parts. After the first step we get that 
 
\begin{align*}
&\mathcal{W}^{v,z}(x)= \sum_{\substack{n \leq x \\ P^{-}(n) >z }} P_v(n)  
=    \sum_{\substack{a \leq x/z \\ P^{-}(a) >z }} P_{v-1}(a) \sum^{*}_{ \substack{ n_1 \leq x/a }} \overline{\theta}(n_1) \\
&=  \sum_{D_1 \in  \mathscr{D}_z(x)} \sum_{\substack{a \leq x/D_1 \\ P^{-}(a) >z }} P_{v-1}(a) \sum^{*}_{ \substack{ n_1 \leq x/a \\ n_1 \sim D_1 }} \overline{\theta}(n_1) \\
& =  \sum_{D_1 \in  \mathscr{D}_z(x)}  \sum_{\substack{a \leq x/D_1 \\ P^{-}(a) >z }} P_{v-1}(a)
\Big[\int_{D_1}^{2D_1} \sum^{*}_{ \substack{ n_1 \leq \min(u_1,x/a) }} \Lambda(n_1) \frac{du_1}{u_1 (\log u_1)^2} \\
&  +\frac{1}{\log (2D_1)} \sum^{*}_{ \substack{ n_1 \leq \min(x/a,2D_1)  }} \Lambda(n_1)  - \frac{1}{\log (D_1)} \sum^{*}_{ \substack{ n_1 \leq \min(x/a,D_1)  }} \Lambda(n_1)  \Big]. 
\end{align*}
Here, the summation $\sum_{n}^{*}$ indicates the restriction $P^{-}(n) >z$.  The above procedure may be repeated inductively,  each time extracting a factor $n_{k+1}$ from the convolution $P_{v-k}(a)$.  Observe that after the $k$-th step,  the dyadic  summation range of $n_{k+1}$ may be chosen in such a way that $D_{k+1}\leq (x/ D_1\cdots D_k)$. In this manner we arrive at a decomposition 
\begin{equation}\label{Wqdecomp}
\mathcal{W}^{v,z}(x)=\sum_{\underline{D} \in \Theta_{v}(x;z) \cap \mathscr{D}_z(x)^v }\  \sum_{ \substack{\wp=\left\{\underline{j}, \underline{k},\underline{l}  \right\} \\ \underline{j}  \dot{\cup}  \underline{k} \dot{\cup} \underline{l} = [v]  } } \int_{ \mathscr{I}_{\underline{j} } (\underline{D}) } \mathcal{W}_{\wp}^{v,z}(x; \underline{u}) \ \beta(\underline{u}) \ d \underline{u},
\end{equation}
which runs over $v$-tuples of dyadic powers $\underline{D}=(D_1,...,D_v)$ in the set 
\begin{equation}
\Theta_{v}(x;z)=\bigg\{\underline{t} \in [z,x]^v: \prod_{j \leq v} t_j \leq 2 x \bigg\}
\end{equation}
and partitions $\wp=\left\{\underline{j}, \underline{k},\underline{l}  \right\}$ such that $ \underline{j} \dot{\cup} \underline{k} \dot{\cup} \underline{l} = [v]  = \left\{ 1,...,v \right\} $. The members of the partition will be denoted by $\underline{j}=\left\{j_1,...,j_r \right\} $,  $\underline{k}=\left\{k_1,...,k_a \right\} $ and $\underline{l}=\left\{l_1,...,l_b \right\} $.  The domain of integration on the RHS of \eqref{Wqdecomp} is 
\begin{equation*}
\mathscr{I}_{\underline{j} } (\underline{D})= \bigtimes_{\ell \leq r} \halfopen{D_{j_{\ell}} }{2 D_{j_{ \ell}} } 
\end{equation*}
and for each tuple of variables $\underline{u}=(u_{j_1},...,  u_{j_r} ) \in \mathscr{I}_{\underline{j} } (\underline{D})$,  the integrand is a product of the functions  
\begin{equation}\label{betadef}
\beta(\underline{u})=\beta(\underline{u}, \wp, \underline{D})= (-1)^a \Big( \prod_{s \leq r} \frac{1}{u_{j_{s}} (\log u_{j_{s}} )^2} \Big) \Big( \prod_{s \leq a} \frac{1}{\log  D
_{k_{s}} } \Big)   \Big( \prod_{s \leq b} \frac{1}{\log  (2D_{l_{s}}) } \Big) 
\end{equation}
and
\begin{equation*}
 \mathcal{W}_{\wp }^{v,z}(x; \underline{u}) =  \sum_{\substack{n_1,...,n_v  \\ n_1 \cdots n_v \leq x} }^{\dagger}  \prod_{s \leq v} \Lambda(n_s).
\end{equation*}
The sum in the last line runs over variables $n_1,...,n_v$ which are subject to the conditions $P^{-}(n_s)>z$ as well as the restrictions 
\[ \forall f \in \underline{j}, \forall t \in \underline{k}, \forall e \in \underline{l}:\  \ \ n_f \leq u_f, \ \ n_t \leq D_t, \ \ n_e \leq 2D_e. \]
Further,  for $x \geq 2y$,  we obtain a corresponding decomposition for $\mathcal{W}^{v,z}(x,y)$ by replacing the long sums $ \mathcal{W}_{ \wp}^{v,z}(x; \underline{u})$ on the RHS of \eqref{Wqdecomp} with their short-interval counterparts 
\[ \mathcal{W}_{\wp }^{v,z}(x,y; \underline{u})= \mathcal{W}_{\wp}^{v,z}(x+y; \underline{u})- \mathcal{W}_{\wp}^{v,z}(x; \underline{u}).\] 

In subsequent sections it will be useful to have a similar expression for the weighted sum $\mathcal{F}_{\underline{\lambda}}(x,y)$. In this case we may simply separate each variable appearing in the convolution sum according to its dyadic range. Proceeding exactly as before (without summing by parts), we get the somewhat simpler identity   
 
\begin{equation}\label{Fdyadicdecomp}
\mathcal{F}_{\underline{\lambda}}(x,y)= \sum_{\underline{D} \in \Theta_{v}(x;z) \cap \mathscr{D}_z(x)^v }\  \sum_{ \substack{\wp=\left\{\underline{k},\underline{l}  \right\} \\  \underline{k} \dot{\cup} \underline{l} = [v]  } }  (-1)^a \ \mathcal{F}_{\underline{\lambda},  \wp}(x,y),
\end{equation}
where for each partition $\wp=\left\{\underline{k},\underline{l}  \right\} $ of the set $[v]$, the sum
\begin{equation*}
 \mathcal{F}_{\underline{\lambda},  \wp }(x,y) = \sum_{m \sim U}  \sum_{\substack{ n_1,...,n_v  \\ m n_1 \cdots n_v \in \halfopen{x}{x+y}} }^{\wedge}  \prod_{s \leq v} \Lambda(n_s)
\end{equation*}
satisfies the restrictions $P^{-}(n_s)>z$ and
\[ \forall t \in \underline{k}, \forall e \in \underline{l}:\   \ n_t \leq D_t, \ \ n_e \leq 2D_e. \]

\subsubsection{Step two: multi-dimensional Perron formulas for the components $\mathcal{W}_{\wp}^{v,z}$ and $ \mathcal{F}_{\underline{\lambda},  \wp }$.}

Suppose that we are given an arbitrary partition $\wp=\left\{\underline{j}, \underline{k},\underline{l}  \right\}$,  together with a $v$-tuple $\underline{D} \in \Theta_{v}(x;z) \cap  \mathscr{D}_z(x)^v $ and a corresponding real tuple $\underline{u} \in  \mathscr{I}_{\underline{j} } (\underline{D})  $. Then we may rewrite the short-interval component $\mathcal{W}_{\wp}^{v,z}(x,y; \underline{u})$ as a $(v+1)$-dimensional Perron formula 
\begin{equation}\label{WcompPerron}
\mathcal{W}_{\wp}^{v,z}(x,y; \underline{u})
= \frac{1}{(2 \pi i)^{v+1}}  \int_{(\kappa)^v}  \int_{(c)} \Xi(\underline{w})   \xi(s,\underline{w}) 
 \mathcal{T}_{x,y}(s)  \frac{ds}{s}   d \underline{w}
\end{equation}
for any $c>1$ and $\kappa>0$. Here the integrand is a product of meromorphic functions defined over the complex variables $s,w_1,...,w_v$. We have 
\begin{equation}\label{Xidef}
\Xi(\underline{w})=\Xi_{\wp, \underline{u},\underline{D}}(\underline{w})=\bigg( \prod_{f \in \underline{j} } \frac{u_{f}^{w_f} }{w_f} \prod_{t \in \underline{k} } \frac{D_t^{w_t} }{w_t}    
\prod_{e \in \underline{l} } \frac{(2D_e)^{w_e} }{w_e}  \bigg),
\end{equation}

\begin{equation*}
\mathcal{T}_{x,y}(s)=(x+y)^s -x^s, \qquad  \ \ Q_{U}(s)=\sum_{ m  \sim U} m^{-s},
\end{equation*}
and
\[\xi(s,\underline{w})=   \prod_{k \leq v} \Big(-\frac{\zeta'}{\zeta}(s+w_k)-\mathcal{L}^{\star}_z(s+w_k) \Big). \]
In the last line we have a shifted product of the restricted log derivative $(-\zeta'/\zeta(s)-\mathcal{L}^{\star}_z(s))$,\footnote{The Dirichlet polynomial $\mathcal{L}^{\star}_z(s)$ is not to be confused with the expression $\mathcal{L}_a(x)$ used elsewhere in this article.} where
\[ \mathcal{L}^{\star}_z(s)=\sum_{k\geq 1, p \leq z} (\log p) p^{-ks}.\] 

The corresponding expression for $ \mathcal{F}_{\underline{\lambda},  \wp }(x,y) $ is given by
\begin{equation}\label{FcompPerron}
 \mathcal{F}_{\underline{\lambda},  \wp }(x,y) 
= \frac{1}{(2 \pi i)^{v+1}}  \int_{(\kappa)^v}  \int_{(c)} \Xi(\underline{w}) Q_U(s)   \xi(s,\underline{w}) 
 \mathcal{T}_{x,y}(s)  \frac{ds}{s}   d \underline{w},
\end{equation}
where, in this setting,  $\underline{j}$ is taken to be the empty set in \eqref{Xidef} and $\wp=\left\{\underline{k},\underline{l}  \right\} $.
\subsection{A discussion of the Dirichlet series}
In this section we collect some useful expansions and estimates for the generating series of $\mathcal{F}_{\underline{\lambda}}$ and  $\mathcal{W}^{v,z}$ in the critical strip.  
\subsubsection{The zero density estimate}
The "$y$ aspect" in the \textsl{comparative asymptotic evaluation} of $\mathcal{F}_{\underline{\lambda}}$ and  $\mathcal{W}^{v,z}$ will depend crucially on the best available zero density estimate for $\zeta$. To state the density result, write
\[ N(\sigma,T) = |\left\{\rho=\beta+i \gamma \in \C: \beta \geq \sigma, |\gamma| \leq T, \ \ \zeta(\rho)=0  \right\} |\] 
for each $\sigma \geq 1/2$ and $T >0$, say. Then we have Guth and Maynard's zero density estimate \cite[Theorem 1.2]{GM}
\begin{equation}\label{Hux}
N(\sigma,T)  \ll T^{\frac{30}{13} (1-\sigma) +o(1)} .  
\end{equation}

\subsubsection{The treatment of $Q_{U}$} Next we recall Vinogradov's estimate for exponential sums with large phase.  This result implies that for values of $s$ with large imaginary part and lying sufficiently close to the $1$-line, the sum $Q_{U}(s)$ has cancellation. 

\begin{lem}[Vinogradov's estimate  \cite{IK}, Corollary 8.26] 
For $t \geq X \geq 2$ we have the bound
\begin{equation}\label{vino}
\big| \sum_{n \leq X} n^{it} \big| \ll X \exp \left( - \kappa \frac{ \log^3 X }{  \log^2 t} \right),  
\end{equation}
where $\kappa>0$ is an absolute constant.  As a consequence, if $M \leq (\log U)^{1/2} $, we have for each $s=\sigma+it$ with imaginary part in the range $ 2 \le U^{1/4} \leq |t| \leq U^{M}$, that
\begin{equation}\label{Qbound}
|Q_{U}(s)|\ll U^{1-\sigma- \kappa'/ M^2}.
\end{equation}

\begin{proof}
The estimate \eqref{Qbound} follows from the claim
\begin{equation}\label{nuvinogradov}
\big|\sum_{  n \in [U, U'] }  n^{it} \big| \ll U  \exp\Big(- \kappa' \frac{\log U}{ M^2 } \Big), \qquad U' \in [U,2U],
\end{equation}
upon summing by parts. To see why \eqref{nuvinogradov} holds,  apply \eqref{vino} in the range $U\leq |t| \leq U^{M}$ and van der Corput's third derivative estimate \cite[Corollary 8.19]{IK} in the range $U^{1/4} \leq |t| \leq U$. 
\end{proof}
\end{lem}

\subsubsection{Some bounds for $\zeta'/\zeta$ in the critical strip.}\label{zetasection}
To estimate the size of $(\zeta'/\zeta)^{v}$ and its variant $\xi(s,\underline{w})$, we start with the well-known approximate formula \cite[Theorem 9.6(A)]{Tit} 
\begin{equation}\label{localzetarep}
\frac{\zeta'}{\zeta}(s)=\sum_{ \substack{ \rho=\beta+i \gamma\\ |\gamma-t|\leq 1 }} \frac{1}{s- \rho} +O(  \log (2+|t|) ), \qquad Re(s) \in [-1,2].
\end{equation}
Since there are at most $O(\log (2+|t|))$ zeroes in the above sum we find that, around each $\rho$ of height $|\gamma| \asymp t$, there exists a radius $( \log (2+|t|))^{-1} \ll C_{\rho} \ll  \log (2+|t|)^{-1}$ such that  
\begin{equation}\label{localzetasup}
\sup_{s:|s-\rho| =C_{\rho}} \Big| \frac{\zeta'}{\zeta}(s) \Big| \ll (\log(2+|t|))^2.
\end{equation}
As an immediate consequence of this local supremal bound we are able to select piecewise linear paths which avoid zeta zeroes in the critical strip.  Such paths will be important later on, in Section 4, when dealing with the Perron formula representations of $\mathcal{F}_{\underline{\lambda}}(x,y)$ and $\mathcal{W}^{v,z}(x,y)$. Before moving on to the description of the paths,  we first take note of the straightforward bound 
\begin{equation*}
|\mathcal{L}^{\star}_z(a)| \ll z^{1-\beta} (\log z), \qquad z \geq 3, \qquad Re(a)=\beta  \in [1/2,1].
\end{equation*}\\
{\bf The contours $\alpha_{\underline{w}}$ and $\Gamma^{\pm}_{\underline{w}}$. } Suppose that we are working with an arbitrary ordinate $\sigma_0 \in [1/2,1]$ and large parameters $x \geq T >0$. Set $\kappa=1/\log x$ and suppose, moreover,  that we have fixed a tuple $\underline{w}=(w_1,...,w_v) $ of complex numbers with each component $w_j$ lying in the narrow rectangle $\left\{ Re(w_j) \in [\kappa, 2\kappa], \ |Im(w_j)| \leq T \right\}$.  Taking advantage of \eqref{localzetasup}, we may select (or rather deduce the existence of) a three-part contour which stays at a distance of at least $\log(2+|T|)^{-2}$ from all shifted zeta zeroes $(\rho - w_j)_{j \leq v, \zeta(\rho)=0}$.  

\noindent The first part, denoted by $\alpha_{\underline{w}}$,  will consist of a piecewise linear path connecting some pair of points $s^{-}= \sigma_0^{-} -i T^{-} $ and $s^{+}= \sigma_0^{+} +i T^{+} $.  Although the precise locations of $s^{-}$ and $s^{+}$ are unimportant, and even impossible to control, we do wish to impose the restrictions $\sigma_0^{\pm} \in [\sigma_0, \sigma_0+ 1/ \log T]$  and $T^{\pm} \in [T/4,4T]$.  The second and third part of the contour are given by the horizontal line segments $\Gamma^{-}_{\underline{w}}=[\sigma_0^{-} -i T^{-}, 1 -i T^{-}  ]$ and $\Gamma^{+}_{\underline{w}}=[\sigma_0^{+} +i T^{+}, 1 +i T^{+}  ]$.  Arguing as in \eqref{localzetasup}, and recalling the pointwise bound for $|\mathcal{L}^{\star}_z(a)|$ stated just below it, we obtain a contour satisfying the following two properties.  First, the path $\alpha_{\underline{w}}$ is made up of at most $O(T(\log T)^2)$ horizontal and vertical line segments and is entirely contained in the strip $\sigma_0 \leq Re(s) \leq \sigma_0 + 1/ \log T$.  Second, we have the estimate
\begin{equation}\label{contoursupboundxi} 
 |\xi(s, \underline{w})| \ll (C \log (2+|t|))^{3v} z^{(1-\beta)v}  (\log z)^v , \qquad \ s=\beta+i t \in \alpha_{\underline{w}} \cup \Gamma^{-}_{\underline{w}} \cup \Gamma^{+}_{\underline{w}}. 
\end{equation}

{\bf The treatment of $\xi(s, \underline{w})$.} Our next task is to obtain level-set estimates for the function  $\xi(s,\underline{w})$.  From \eqref{localzetarep},  we obtain the estimate (outside a null set of values of $\underline{w}$)
\begin{equation}\label{xibounddistance}
|\xi(s,\underline{w})| \leq z^{v(1- \beta)} \prod_{j\leq v}   \Big[B' \log(2+z+|t|) \Big(1+ \frac{1}{\mathfrak{m}_{\zeta}(s+w_j) } \Big) \Big],
\end{equation}
where $B'>0$ is some absolute constant, $Re(s)=\beta$ and $\mathfrak{m}_{\zeta}(a)$ denotes the distance from the complex number $a$ to the nearest zeta zero.  In the following lemma, we will show that $|\xi(s,\underline{w})|$ is not too large \textsl{on average} and exploit this fact to obtain some residue estimates. 

\begin{lem}\label{avreslemma}
Let $x \geq y \geq 1$ be sufficiently large and set $\kappa=1/ \log x$. Suppose, moreover, that $z \geq 1$ and $v \leq \log x$ satisfy $z^v \leq x$. Then for any $1 \leq T \leq x$ there exists a sequence of ordinates $(\kappa_1, ..., \kappa_v) \in [\kappa, 2 \kappa]^v$ such that
\begin{equation}\label{Taylorrhosans} 
  \int_{ (\underline{\kappa})_T }  \Big| \sum_{l \leq v} \mathop{Res}_{s = \rho -w_l} \Big(    \xi(s,\underline{w}) 
  \frac{\mathcal{T}_{x,y}(s) }{s}  \Big) \Big|  \frac{d \underline{w} }{w_1 \cdots w_v} \ll y (x/z^v)^{\beta-1} (B \log x)^{4v}
\end{equation} 
and
\begin{align}\label{Taylorrho} 
\begin{split}
 \int_{(\kappa_1)_T} & \Big|  \int_{ \widehat{(\underline{\kappa})_T } }  \mathop{Res}_{s = \rho -w_1} \Big(  Q_{U}(s)    \xi(s,\underline{w}) 
  \frac{\mathcal{T}_{x,y}(s) }{s}  \Big) \frac{d w_2 \cdots d w_v }{w_2 \cdots w_v}  \Big|  \frac{d w_1}{w_1} \\
& \ll y (x/z^v)^{\beta-1} (B \log x)^{4v} |Q_U(\rho -w_1)|
\end{split}
\end{align} 
for some absolute constant $B>0$ and any zeta zero $\rho= \beta + i \gamma$ lying in the rectangle $R_T=\left\{s: Re(s) \in [0,1],  |Im(s)| \leq T \right\}$.  Here we have written $(\kappa_l)_T=[\kappa_l-iT, \kappa_l + iT] $ and 
\[ (\underline{\kappa})_T=\bigtimes_{r \leq v} (\kappa_r)_T, \qquad \  \ \widehat{(\underline{\kappa})_T }= \bigtimes_{2 \leq r \leq v} (\kappa_r)_T. \]
\begin{proof}
Let us first consider \eqref{Taylorrhosans}. Since there are at most $O(x \log x)$ zeta zeroes in $R_T$, we first deduce the existence of parameters $(\kappa_1, ..., \kappa_v) \in [\kappa, 2 \kappa]^v$ satisfying 
\[ \sup_{s \in R_T} \sup_{l \leq v} \sup_{w_l \in (\kappa_l)_T} \ |\mathfrak{m}_{\zeta}(s+w_l)^{-1}| \leq 2 x \log x. \]
It will be enough to estimate for each $s \in R_T$, the logarithmic measure of the level set

\[  E_{\underline{K} }(s)= \left\{ (w_2,...,w_v) \in  (\kappa_2)_T \times\cdots \times  (\kappa_v)_T: \forall l\geq 2,  |\mathfrak{m}_{\zeta}(s+w_l)|^{-1} \in[K_l,2K_l] \right\}
\]
for any $v-1$ tuple $\underline{K}=(K_2,...,K_v)$ made up of dyadic components $1 \leq K_l \leq x$.  We begin by observing the trivial bound $ |\mathcal{T}_{x,y}(s) / s| \ll y x^{Re(s)-1}$, which holds for any $s \in R_T$.  At height $t \leq x$, any  unit-length square in the critical strip contains at most $B'' \log x$ many zeta zeroes and thus, 
\[  \int_{ E_{\underline{K} }(s) }   \frac{d w_2\cdots dw_v}{|w_2 \cdots w_v|} \leq \prod_{2 \leq l \leq v} \Big( \frac{B'' (\log x)^2}{K_l} \Big) \]

Since $\xi(s,\underline{w})$ has a simple pole at each shifted zero $\rho - w_1$,  we may now apply \eqref{xibounddistance},  to find that the LHS of \eqref{Taylorrhosans} is at most

\begin{align*}
& \ll v \int_{ (\kappa_1)_T}  \sum_{\underline{K}} \  \sup_{ (w_2,...,w_v) \in  E_{\underline{K} }(\rho - w_1) }  |\xi(s,\underline{w})|  \prod_{2 \leq l \leq v} \Big( \frac{B'' (\log x)^2}{K_l} \Big) \frac{d w_1}{w_1} \\
& \ll y (x/z^v)^{\beta-1}  \sum_{\underline{K}}  (3B'B'' \log x)^{3v} \leq y (x/z^v)^{\beta-1} (B \log x)^{4v}. 
\end{align*}
with $B=3B'B''$.

The proof of \eqref{Taylorrho} follows the exact same argument.
\end{proof}
\end{lem}

{\bf The pole at $\rho=1$.} To conclude this section, we briefly consider the coefficients in the Laurent expansion of $(-\zeta' / \zeta(s))^v$ at the pole $s=1$.  Cauchy's formulas give rise to the expansion
\begin{equation*}
-\frac{\zeta'}{\zeta}(s)=\sum_{j\geq -1} c_j (s-1)^j, \qquad |c_j|\ll C^{j+1},
\end{equation*} 
with some absolute constant $C>0$. Taking the $v$-th power of this local representation we then arrive at the expression
\begin{equation}\label{logexpansion}
\left(- \frac{\zeta'}{\zeta}(s) \right)^{v}=\sum_{j\geq -v} d_j^{[v]} (s-1)^j, \qquad \  \forall |j|\leq v: |d_j^{[v]}| \ll (C')^v.
\end{equation}
Indeed,  we have that
\begin{align*}
|d_j^{[v]}|& \ll  \sum_{\substack{  l_1,\dots,l_{v}\geq -1 \\ l_1+\cdots+l_{v}= j }} \ \prod_{r \leq v } |c_{l_r}|
  \ll   \sum_{\substack{  l_1,\dots,l_{v}\geq 0  \\ l_1+\cdots+l_{v}= j+v }}   C^{j+v}  \\
&=   \binom{j+2 v -1}{j+v} C^{j+v} \leq (1+C)^{3v}.
\end{align*}
We may use this observation to calculate certain residues arising in the study of $\mathcal{F}_{\underline{\lambda}}$. More precisely,  we wish to calculate the "main term" in the residue
\[R_v(x)= \mathop{Res}_{s=1} \Big[ \left( -\frac{\zeta'}{\zeta}(s) \right)^{v} \frac{x^s}{s} \Big].  \]
Let us write
\[ \frac{x^s}{s}=x \sum_{\ell \geq 0} M_{\ell}(\log x) (s-1)^{\ell}, \qquad \  M_{\ell}(\log x)= \sum_{ \substack{ j,r \geq 0\\ j+r= \ell}} \frac{(\log x)^j}{j!} \left(-1 \right)^{r}.\]
In combination with \eqref{logexpansion}, we gather that
\begin{equation*}
R_v(x)=x \sum_{- v \leq j \leq -1}    d_j^{[v]}   M_{-1-j}(\log x).
\end{equation*}
When $v \leq \log x/ (\log_2 x)$,  the $j=-v$ term dominates the sum and thus
 \begin{equation}\label{Residuevbound}
|R_v(x)| \ll x \frac{(\log x)^{v-1}  }{(v-1)!}.
\end{equation}

\section{Comparative asymptotics for  \texorpdfstring{$\mathcal{F}_{\underline{\lambda}}$ and  $\mathcal{W}$.}{F{lambda}  and V{lambda} } }

\noindent The asymptotics and estimates given in this section form the backbone of our entire strategy.  The following proposition, and its corollary,  allow us to extend the region of validity of \eqref{curlyFmeanmoderateshort} to even shorter intervals and furnish comparative asymptotic formulas for both $\mathcal{F}_{\underline{\lambda}}(x,y)$ and $\mathcal{W}^{v,z}(x,y)$.  

\begin{prop}\label{extensionprop}
a) Let $\varepsilon>0$ and $a>4$ be arbitrary but fixed. Let $x \geq 1$ be sufficiently large and assume that $1 \leq v \leq  \mathcal{L}_a(x) $. Suppose, moreover, that 
\[  \frac{\log x}{ \exp(3 L_v(x)^{9/10} ) \log v} \leq \log  U\leq  \frac{\log x}{\log_3 x}, \qquad \ \  1 \leq z \leq U^{1/3v}.\] 
Then the formula \eqref{curlyFmeanmoderateshort} remains valid in the extended range $x^{17/30 + \varepsilon} \leq y \leq x$.  More precisely, we have that
\begin{equation}\label{curlyFextendedrange}
\mathcal{F}_{\underline{\lambda}}(x,y) = y\sum_{\substack{ l \in  \halfopen{ \frac{x+y}{2U} }{  \frac{x}{U} } \\ P^{-}(l) >z }} \frac{ F_{v}(l) }{l} + O_{\varepsilon }(y).
\end{equation}

b) For the smaller values $1 \leq v \leq (\log x)^{1/3}/(\log_2 x)^2 $, we have the asymptotic formula

\begin{equation}\label{shortWasymp}
\mathcal{W}^{v,z}(x,y) = \frac{y}{x} \mathcal{W}^{v,z}(x) + O_{\varepsilon }\Big(\frac{y}{ (\log x)^{2v}} \Big),
\end{equation}
in addition to the estimate
\begin{equation}\label{shortnuvanilla}
\sum_{x \leq n \leq x+y} F_v(n) \ll y \frac{(\log x)^{v-1} }{(v-1)!}.
\end{equation}

\begin{proof}
{\bf Part $a$.} Let us first examine the bilinear sum $\mathcal{F}_{\underline{\lambda}}(x,y)$.  Set $h=x/y$ where, throughout the entire argument, we assume that $x^{17/30 +\varepsilon} \leq y \leq x(\log x)^{5v}/U$.  Given an arbitrary partition $\wp=\left\{\underline{k},\underline{l}  \right\} $ of $[v]$ together with a $v$-tuple $\underline{D} \in \Theta_{v}(x;z) \cap  \mathscr{D}_z(x)^v $,  it will be enough to furnish a comparative asymptotic for the component $ \mathcal{F}_{\underline{\lambda}, \wp}(x,y) $, by way of the Perron formula \eqref{FcompPerron}.  

Set $\kappa=1/ \log x$ and $T= x^{\varepsilon/ 2} h$,  and let $(\kappa_1, ..., \kappa_v) \in [\kappa, 2 \kappa]^v$ be the sequence of ordinates given in Lemma \ref{avreslemma}.  We may truncate each of the integrals over the $w_j$ variables at height $\pm T$ and then truncate the integral over $s$ at suitable heights $T^{-}_{\underline{w}} ,T^{+}_{\underline{w}} \asymp T$ below and above the real line. In view of the pointwise bound $\sum_{d|n} F_v(d) \ll_{\varepsilon} x^{\varepsilon /4}$ (valid for all $n \leq 2x$),  it follows that
\begin{equation*} 
 \mathcal{F}_{\underline{\lambda}, \wp}(x,y)  
=  \frac{1}{(2 \pi i)^{v+1}}  \int_{(\underline{\kappa})_T}  \int_{1-iT_{\underline{w}}^{-} }^{1+iT_{\underline{w}}^{+} }  \Xi(\underline{w})   \xi(s,\underline{w}) 
Q_U(s) \mathcal{T}_{x,y}(s)  \frac{ds}{s}   d \underline{w}+ O_{\varepsilon}(x^{1+\varepsilon/3}/T).
\end{equation*} 
We now set $\sigma_0=17/30$ and first consider the  above expression for fixed $\underline{w}$, deforming the $s$-contour into three parts, along the lines of the discussion in section \ref{zetasection}. The first part is a piecewise linear path $\alpha_{\underline{w}}$ which is contained in the narrow strip  $Re(s) \in [\sigma_0,\sigma_0+1/\log x]$ and connects the points $\sigma_0^{-}-iT_{\underline{w}}^{-}$ and $\sigma_0^{+}+ iT_{\underline{w}}^{+}$. The two remaining parts are the horizontal line segments $\Gamma^{+}_{\underline{w}}=[\sigma_0^{+}+iT_{\underline{w}}^{+},1+iT_{\underline{w}}^{+}],\ \Gamma^{-}_{\underline{w}}=[1- iT_{\underline{w}}^{-},\sigma_0^{-}-iT_{\underline{w}}^{-}]$.  

To carry out the integration over $\alpha_{\underline{w}} \cup \Gamma^{-}_{\underline{w}}\cup \Gamma^{+}_{\underline{w}}$, we recall the estimate \eqref{contoursupboundxi} and the trivial bounds $|\Xi(\underline{w})|\leq 2^v/|w_1\cdots w_v| $ and $|Q(s)|\leq U^{1-\sigma} $ to obtain a contribution no greater than

\begin{align*}
& \ll \int_{(\underline{\kappa})_T} (C \log x)^{4v} \Big[ x^{\sigma_0} (z^v U)^{1-\sigma_0} \int_{\alpha_{\underline{w}}}    \frac{ \ d s }{|s|} +  \Big(  \sup_{\beta \in [\sigma_0, 1]} \frac{ x^{\beta} U^{1-\beta} z^{v (1- \beta)}  }{T} \Big) \Big] \frac{d \underline{w}}{|w_1\cdots w_v|} \\
& \ll (x^{\sigma_0} U z^{v(1-\sigma_0)} + x/T) (C \log x )^{5v +1} =O(yx^{-\varepsilon /4}).
\end{align*}

As a result of the deformation we also pick up the residues 
\begin{align*}
\Upsilon(x,y)&= \frac{1}{(2 \pi i)^{v}}  \sum_{j \leq v} \int_{(\underline{\kappa})_T}  \sum_{\substack { \rho   \\   \rho -w_j \in \mathcal{D}_{\underline{w} }  } }   \Xi(\underline{w}) \mathop{Res}_{s = \rho -w_j} \left(    \xi(s,\underline{w}) Q_U(s)    \frac{\mathcal{T}_{x,y}(s) }{s}\right) d \underline{w} \\
&=  \frac{1}{(2 \pi i)^{v}}  \sum_{j \leq v}  \big(\Upsilon^{-}_j(x,y) +\Upsilon^{+}_j(x,y) \big).
\end{align*}
Here,  $\rho$ runs over zeta zeroes and we have used the disjointedness of the sets $Z_j=\left\{ \rho -w_j: \zeta(\rho)=0 \right\}$ to separate the residues.  We have written $\mathcal{D}_{\underline{w}}$ to denote the subregion of the critical strip enclosed by the $\alpha_{\underline{w}} \cup \Gamma^{-}_{\underline{w}} \cup \Gamma^{+}_{\underline{w}}$ and the line $Re(s)=1$. The sum  $\Upsilon^{+}_j(x,y)$ denotes the sum over residues restricted to shifted zeta zeroes $\rho-w_j=\beta+i \gamma$ in the range $ U^{1/3} \leq |\gamma| \leq T$ and $\Upsilon^{-}_j$ counts all the remaining residues in $w_j +\mathcal{D}_{\underline{w}}$. 
\\

{\bf The estimation of $\Upsilon_{j}^{+}(x,y)$. }To treat the high-lying zeroes $\rho-w_j=\beta+i \gamma$ satisfying $U^{1/3} \leq |\gamma| \leq T= x^{\varepsilon/2} h $, we set $M=\exp(3 L_v(x)^{9/10} ) \log v$ and appeal to \eqref{Hux} and \eqref{Qbound}.
For zeroes with large enough real part, that is $\beta \geq \sigma':=1- \kappa'/2M^2$,  we have that 
\[ |Q_{U}(\rho)|\ll U^{1-\beta- \kappa'/M^2} \ll U^{- \kappa'/2M^2}. \]
Combined with the zero-density estimate \eqref{Hux}, and the assumption that $z^{v} \le U^{1/3}$, we may now invoke the bound \eqref{Taylorrho} to find that 
\begin{align*}
& |\Upsilon^{+}_j(x,y)|  \ll (\log x)^{4v} y  \Big( \sum_{\substack{ \rho-w_j:  \ \beta \in [\sigma_0, \sigma']\\  \ U^{1/3} \leq |\gamma|< T \\ }} (x/Uz^{v})^{\beta-1}
+ \sum_{\substack{ \rho - w_j: \ \beta \in [\sigma',1]\\  \ U^{1/3} \leq |\gamma|< T }}  (x/z^{v})^{\beta-1} U^{- \kappa'/2M^2} \Big) \\
& \leq (\log x)^{4v} y \left( \int_{\sigma_0}^{\sigma'} (x/ Uz^{v})^{\beta -1} \ d N(\beta,T)  +   \int_{\sigma'}^{1}   (x/z^{v})^{\beta-1}  U^{- \kappa'/2M^2} \ d N(\beta,T)  \right) \\
& \ll (\log x)^{4v+1} y \Big[  \sup_{\beta \in [\sigma_0, \sigma']} \big( \frac{x}{Uz^{v}} \big)^{\beta -1} N(\beta,T)  +   \sup_{\beta \in [\sigma',1 ]}   \frac{(x/z^{v})^{\beta-1}  N(\beta,T)  }{U^{\kappa'/2M^2} }  \Big] =O\Big( \frac{y}{(\log x)^{2v} } \Big).
\end{align*}
In the last line we used the estimate $(\log x)^{7v} U^{- \kappa'/2M^2} =O(1) $ together with the bound $T^{30/13} =o(x)$. To see why the former estimate holds, it is enough to verify that $M^2 \leq (\kappa'/14)(\log U / v \log_2 x) $, which in turn follows from $M^3 \leq (\kappa'/14)(\log x / v \log_2 x) $. To prove this last inequality we observe that
\[ \log M=\log_2 v+ 3 L_v(x)^{9/10} <  \frac{1}{3} \big[ L_v(x) + \log \frac{\log (v+1)}{\log_2 x} \big] \]
for all large $x \geq 1$ and $v \leq \log x/ (\log_2 x)^a$ (recall that $a>4$).\\

{ \bf The treatment of $\Upsilon^{-}_j(x,y)$. } In this last part we set up a comparison of the quantities $ \mathcal{F}_{\underline{\lambda}}(x,y) $ and $ \mathcal{F}_{\underline{\lambda}}(x,y^{*}) $ where $y^{*}=x(\log x)^{5v} /U$.  Let us restrict our attention to $\Upsilon^{-}_1(x,y)$ for notational convenience.  Since $\xi(s, \underline{w})$ has a simple pole at each shifted zero $\rho_1=\rho - w_1$, we first observe that
\[ \mathop{Res}_{s = \rho_1} \left(     \xi(s,\underline{w}) Q_U(s)    \frac{\mathcal{T}_{x,y}(s) }{s} \right)= \xi(\rho_1, \widehat{\underline{w} }) Q_U(\rho_1)    \frac{\mathcal{T}_{x,y}(\rho_1) }{\rho_1},  \]  
where $ \widehat{\underline{w} }=(w_2,...,w_v)$. The  desired comparison follows upon writing a Taylor expansion in the form
\begin{align*}  \mathcal{T}_{x,y}(\rho_1)&=\frac{(x+y)^{\rho_1}  -  x^{\rho_1}  }{\rho_1}= yx^{\rho_1 -1} +\frac{y^2}{2}(\rho_1 -1) x^{\rho_1 -2}+\cdots\\
&=: y B(x, \rho_1) + R(x,y,\rho_1). 
\end{align*} 
Here the expression $yB(x, \rho_1)$ corresponds to the linear term in the expansion. The remainder term satisfies the bound
\[|R(x,y,\rho_1)| \ll y x^{\beta -1} \frac{|\rho_1|}{h} \ll  \frac{ y x^{\beta -1} }{U^{1/2}}   , \qquad |\gamma|< U^{1/3},\] 
where,  in the last line, we have used that $(\log x)^v=U^{o(1)}$.  By (the proof of) \eqref{Taylorrho} we have that 
\begin{align*}
 \Upsilon^{-}_j(x,y)& = \int_{(\underline{\kappa})_T }  \Xi(\underline{w}) \sum_{\substack{ \rho -w_1 \in D_{\underline{w}}\\  \ \beta \geq \sigma_0 \\  \ |\gamma|< U^{1/3}  \\ }}   \Big(    \xi(\rho_1, \widehat{\underline{w} }) Q_U(\rho_1)  
   \left[  y B(x, \rho_1)+ O\big( \frac{ y x^{\beta -1} }{U^{1/2}} \big)   \right]  \Big) d \underline{w} \\
& =:  y G(x)+O \left(    \frac{(\log x)^{4v} y }{U^{1/2}}   \int_{\sigma_0}^{1}   (x/Uz^v)^{\beta -1}  \ d N(\beta,T)    \right). 
\end{align*}
Using once again that $(\log x)^v = U^{o(1)}$ we see,  as in the estimation of $ \Upsilon^{+}(x,y)$,  that the error term in the last line is $O(y/(\log x)^{2v})$ and hence, gathering all of our estimates so far, we find that 
\[  \mathcal{F}_{\underline{\lambda}, \wp }(x,y) = y G(x)+ O\Big( \frac{y}{(\log x)^{2v} } \Big)= \frac{y}{y^{*}}\mathcal{F}_{\underline{\lambda}, \wp}(x,y^{*}) + O\Big( \frac{y}{(\log x)^{2v} } \Big), \]
uniformly over all $x^{17/30 +\varepsilon} \leq y \leq x(\log x)^{5v}/U$, which concludes the treatment of $\mathcal{F}_{\underline{\lambda}, \wp}(x,y)$.  

All that remains is to sum the comparative statement,  given just above, over all partitions $\wp=\left\{\underline{k},\underline{l}  \right\} $ of $[v]$ and all $v$-tuples $\underline{D} \in \Theta_{v}(x;z) \cap  \mathscr{D}_z(x)^v $.  In doing so, we multiply the error by a factor no greater than $(2 \log x)^v$ and thus we retrieve \eqref{curlyFextendedrange}.  \\ \\

\noindent {\bf Part $b$.} The proof of \eqref{shortWasymp} is very similar to the discussion in part a, and hence we offer the following summary, highlighting the key facts.  In the present setting we make use of \eqref{WcompPerron} in place of \eqref{FcompPerron} and proceed, as before, to deform the $s$-contour for each fixed tuple $\underline{w}$.  After estimating the integral along $\alpha_{\underline{w}} \cup \Gamma^{-}_{\underline{w}}\cup \Gamma^{+}_{\underline{w}}$ (which yields the exact same bound $O(yx^{-\varepsilon /4})$) we are left with the residue sum $\Upsilon(x,y)$.  Invoking the Vinogradov-Korobov zero-free region, we may separate the pole at $\rho =1$ from those in the region 
\[R_{VK}(T)=\left\{ \beta \leq 1- \frac{1}{(\log T)^{2/3} (\log_2 T)^{3/4} }= :\overline{\sigma} \right\}. \] 
Correspondingly,  we denote by $\Upsilon_0$ the restriction of $\Upsilon(x,y)$ to poles lying in the region $R_{VK}(T)$.  As in the estimation of $|\Upsilon^{+}_j(x,y)| $ in part a,  we find that
\begin{equation}\label{RVKbound}
|\Upsilon_0| \ll y (\log x)^{4v+1}  \sup_{\beta \in [\sigma_0, \overline{\sigma}]}  (x/z^v)^{\beta -1} N(\beta,T) \ll y (\log x)^{4v} \Big( \frac{T}{x^{13/30}} + x^{\overline{\sigma}-1 }  \Big), 
\end{equation} 
incurring an error of size at most $O(y(\log x)^{4v+1} x^{\overline{\sigma}-1 } )=O(y/ (\log x)^{2v})$. Finally, the pole at $s=1$ is dealt with as in the treatment of $ \Upsilon^{-}_j(x,y)$, yielding \eqref{shortWasymp}.\\

To get \eqref{shortnuvanilla} first observe that, thanks to the above argument,  we also obtain the comparative asymptotic equivalence 
\[ \sum_{x \leq n \leq x+y} F_v(n) \sim \frac{y}{x} \sum_{n \leq x} F_v(n).  \]
Since
\[ \sum_{ n \leq x} F_v(n) = \frac{1}{2 \pi i} \int_{(c)}  \Big(-\frac{\zeta'}{\zeta}(s) \Big)^{v} x^s \frac{ds}{s} + O_{\varepsilon}(x^{1+\varepsilon/3}/T), \]
we may proceed,  as before, to deform the contour along the boundary $\partial(R_{VK}(T) )$, where the integrand satisfies the bound $O(x/|s|)$. The main contribution, appearing on the RHS of \eqref{shortnuvanilla} comes from the residue at  $s=1$, which may be bounded by \eqref{Residuevbound}.
\end{proof}

\end{prop}

As a consequence of the proposition we are now in a position to produce a sufficiently sharp estimate for $\mathcal{F}_{\underline{\lambda}}(x,y)$.  One last ingredient is required.  

\begin{lem}\label{combilemma}
For $v \geq 1$, set $K=\lfloor \log v / \log 2 \rfloor$ and define the combinatorial product 
\begin{equation}\label{combiJdef}
 \mathcal{C}_{v}(K)= \max_{w \leq v} \max_{ \substack{ a_0,...,a_K \\ w+ a_0+...+a_K \leq v \\  \forall j \geq 0: a_j  \leq v 2^{-j}  }  } { {v }\choose{ w} }  \prod_{0 \leq k \leq K} { {v }\choose{ a_k} } .
\end{equation}
Then we have the inequality $ \mathcal{C}_{v}(K) \leq A^v$, where $A=16\sqrt{2e}$ is admissible.
\begin{proof}
We make use of the standard estimates $\max_w { {n }\choose{ w} }  \leq 2^n $  and ${ {n }\choose{ r} } \leq (en/r)^r $ to find that
 \[ \mathcal{C}_{v}(K) \leq \max_{w \leq v} { {v }\choose{ w} }  \prod_{0 \leq k \leq K} \left( \max_{a_k \in [0,\frac{v}{2^{k}}]} { {v }\choose{ a_k} } \right) \leq 8^{v} \prod_{2 \leq k \leq K} \Big( \frac{e v}{v/2^k} \Big)^{v/2^k} = (16\sqrt{2e})^{ v}.
\]

\end{proof} 
\end{lem}

\begin{cor}\label{Evalcor}
Given sufficiently large parameters $x \geq 1$ and $\log_3 x \leq v \leq \mathcal{L}_2(x) $, we have the estimate
\begin{equation}\label{Fnumean}
 \sum_{\substack{  n \leq x\\ P^{-}(n)>z  }} \frac{ F_{v}(n) }{n} \ll  v \left( \frac{\kappa \log x }{v}  \right)^{v}, 
\end{equation}
provided that  $z \geq 30 v $. Here $\kappa = 150$ is an admissible constant. 
Moreover,  under the assumptions of part a of Proposition \ref{extensionprop}, we have that
\begin{equation}\label{curlyFmean}
\mathcal{F}_{\underline{\lambda}}(x,y) \ll y \log x \left( \frac{ \kappa \log x }{v } \  \right)^{v}. 
\end{equation}
\begin{proof}
Let us write $v_x'= (\log x)/v$, as before, and begin with the treatment of \eqref{Fnumean}. In view of \eqref{Fvsupport}, it suffices to fix a $ w \leq v$ and obtain estimates for the sum \eqref{Fnumean} restricted to integers $n \in S_w$. To this end, we separate each $n \in S_{w}(x)$, into two parts:
\[n_0= \prod_{p^e ||n, \  p \leq \exp(v_x')} p^e, \ \qquad n_1=n/n_0,  \]
setting $w_0=\omega(n_0)$ and $w_1=\omega(n_1)$ with $w_0 +w_1=w$. 
In order to exert more control over the size of $F_{v}(n)$, we will take into account the prime decomposition of $n_1$: write $w_1= \beta v $ and consider the long intervals $J^{-}=[z,\exp(v_x')]$ and $J_k=[\exp(v_x' 2^k), \exp( v_x' 2^{k+1})]$ for $0 \leq k \leq \lfloor \log v / \log 2 \rfloor=:K$. Observe that in each range $J_k$, there are at most $ v 2^{-k}$ many primes dividing $n_1$, so we may work under the assumption that $\omega_{J_k}(n)=\alpha_k v 2^{-k}$ for some fixed sequence $\underline{\alpha}=(\alpha_0,..., \alpha_K) \in [0,1]^{K+1}$. Moreover, there are at most $w^{\log v / \log 2} $ possible choices for $\underline{\alpha} \in [0,1]^{K+1}$. Using the shorthand $\| \underline{\alpha} \|_1^{*}=\sum_{k \geq 0}  \alpha_k 2^{-k}$ and recalling the assumption that $z\geq 30 v$ as well as the combinatorial factor \eqref{combiJdef}, we get the estimate
\begin{align*}
 &  \qquad \sum_{\substack{  n \in S_{w}(x) \\ P^{-}(n)>z  } } \frac{ F_{v}(n) }{n}    
\leq w^{\log v / \log 2} \ \mathcal{C}_{v}(K)  \sup_{ \substack{ \underline{\alpha} \in [0,1]^{K+1} \\ \| \underline{\alpha} \|_1^{*}=\beta } }   \Big\{ \\
  &  \Big( \sum_{p \in J^{-}} \big[ \frac{\log p}{p} +\sum_{r \geq 2}  \frac{ (\log p) v^r}{p^{r} } \big] \Big)^{w_0}  \prod_{ k \geq 0}  \Big( \sum_{p \in J_k} \big[ \frac{\log p}{p} +\sum_{r\geq 2}  \frac{ (\log p) v^r}{p^{r} } \big] \Big)^{\alpha_ k v 2^{-k}}     \Big\} \\
& \leq  v^{2 \log v}A^{v} \sup_{ \substack{ \underline{\alpha} \\ \| \underline{\alpha} \|_1^{*}=\beta } }  \Big\{ \big( v_x' + \frac{v_x' (v/z)^2}{1-v/z} +2 \big)^{w_0} \prod_{ k \geq 0}
\Big( 2^{k} v_x' + \frac{2^{k+1}  v_x' (v/z)^2}{1-v/z} +2\Big)^{\frac{\alpha_ k v}{2^{k}} }\Big\} \\
& \leq v^{2 \log v} A^{v} \sup_{ \substack{ \underline{\alpha} \\ \| \underline{\alpha} \|_1^{*}=\beta } }  \Big\{ \Big( \frac{201}{200} v_x' \Big)^{w_0+ v \| \underline{\alpha} \|_1^{*}} \exp \Big( v \log 2 \sum_{k \geq 1} \frac{k \alpha_ k }{2^{k}} 
 \Big) \Big\} \leq \big( 150  v_x' \big)^{v}, 
\end{align*} 
which retrieves the desired bound \eqref{Fnumean}, upon adding the contribution of each $w \leq v$. The expression $\mathcal{C}_{v}(K)$ bounds the number of appearances of each squarefree divisor of $n$ in the convolution expansion $F_{v}(n)$ - the factors $v^r$ in the second line play a similar role for prime-power divisors of $n$. In the third line we used Mertens' bound $\sum_{p \leq u} \log p/ p \leq \log u +2$.
\end{proof}
\end{cor}

\section{Relating \texorpdfstring{$\pi_{\nu}(x,y)$}{pi{nu}(x,y) } to the short minorant  \texorpdfstring{$ \mathcal{M}^{\sharp}_{\nu} (x,y) $}{M{nu}(x,y)} }\label{upperboundproofsection}

This section is devoted to the proof of Theorem \ref{squarefreemainprop} which asserts that the minorant
\begin{equation}\label{shortminorant}
 \mathcal{M}^{\sharp}_{\nu} (x,y)=\sum_{1 \leq w \leq \ell_{\nu}(x) } \ 
\sum_{ \substack{m \in S_{\nu -w}(\tau) \\ m | \mathcal{Q}_{\nu}(x)^{\infty}} }\ 
\sum_{\substack{ n  \in S_w( x/m,y/m) \\  (n, \mathcal{Q}_{\nu}(x) )=1 } } 1.
\end{equation} 
captures almost all integers in $S_{\nu}(x,y)$. As before we have used the notation $\tau=\exp(\lambda^{+})$.

The proof of Theorem \ref{squarefreemainprop} will be carried out by way of an inductive mechanism.  Let us begin with an outline of the argument.\\
We may assume that $\nu \geq \log_2 x/ (\log_3 x)^{3/4}$-- otherwise we have \eqref{simpleanatomy}. The key procedure in our dissection of $S_{\nu}(x,y)$ is the removal of the \textsl{bad} integers in $S_{\nu}(x,y)$, that is to say integers which do not appear in $\mathcal{M}^{\sharp}_{\nu}(x,y)$.  In order to introduce the bad sets, we require suitable cut-off/ threshold parameters: given $v \leq \nu$ and  $c \in (0,99/100)$, write
\[ z^{*}_v(c,x)= \exp \Big(  \frac{ v'_x}{2 \exp(L_{v}(x)^c ) } \Big). \] 
For a typical $n \in S_{v}(x)$, we expect the divisor 
\[ n_0=\prod_{p^e || n, p \leq z_v^{*}(c,x) } p^e \]  
to contain a large proportion of the prime factors of $n$. Accordingly, when dealing with atypical integers, we define for any $\eta \in (0,1)$ and $t>1$, a (generalised) set of bad integers for which this is not true:
\begin{equation}\label{badsetnudefinition}
\mathcal{B}_{v}(x,y;\eta, c,t)=\left\{ n \in S_{v}(x,y) : \omega(n_0)<(1-\eta) v, \  P^{-}(n)>t  \right\}. 
\end{equation}
In the above notation we thus arrive at a partition 
\begin{equation}\label{goodvsbad} S_{\nu}(x,y) =   \mathcal{G}_{\nu}(x,y;\eta, c,1) \cup \mathcal{B}_{\nu }(x,y;\eta, c,1) \end{equation} 
consisting of a \textsl{good and bad subset} of $S_{\nu}(x,y)$ and seek to estimate the cardinality of the bad set.

\subsection{Estimating the cardinality of the bad set}

\begin{prop}\label{dealingwiththebadset} Let $\eta \in (0,1)$,  and $\varepsilon>0$ be given. Under the same assumptions as Theorem \ref{squarefreemainprop}, we have for any $\log_3 x \leq v \leq \mathcal{L}_a(x)$,  any $c$ as in \eqref{widecrange} and any $1 \leq t \leq (z_v^{*}(c,x))^{1/3}$ the estimate
\begin{equation} \label{badsetstarbound}
|\mathcal{B}_{v}(x,y;\eta,c,t)|  \ll_{\eta,\varepsilon}  y a_{v}(x) K_{\varepsilon}^{v} \Delta^{*}_{v,\eta v/3,x'}(x; c,t) +E_{v}(x,y) 
\end{equation} 
and
\begin{equation} \label{secondbadsetstarbound}
|\mathcal{B}_{v}(x,y;\eta,c,t)|  \ll_{\eta,\varepsilon}  y a_v(x) K_{\varepsilon}^{v} \frac{\delta_v(x)}{L_v(x)^{(1-c) \eta v/ 3} }.
\end{equation}
Here,  any choice of $x' \in [x,x^{4/3}]$ is admissible, 
\[ E_{v}(x,y) \ll \sqrt{x}+ \frac{y}{2^v} \underline{\delta}_{v} (x;t)   \]
and 
\begin{equation*}
a_{v}(x) =
\left\{
	\begin{array}{lll}
		(\log x)^6 & \mbox{if }  v>(\log_2 x)^{10} \\
		& \\
		1 & \mbox{if } v \leq (\log_2 x)^{10}.
	\end{array}
\right.
\end{equation*}
$K_{\varepsilon}>10$ depends only on $\varepsilon$, $\Delta_{v,u,x'}^{*}(x;c,t)$ was defined in \eqref{Deltastardef} and $ \underline{\delta}_{v} $ in \eqref{lowerdeltatdef}.

 \begin{proof}
Throughout the proof of the proposition we will use the notation
\[ x_0=\exp \Big( \frac{\log x}{2 \exp( L_{v}(x)^c ) }  \Big)  \] 
and assume that $x^{17/30 + \varepsilon} \leq y \leq x$. In order to take full advantage of Proposition \ref{extensionprop}, we must first reduce the estimation of $|\mathcal{B}_{v}(x,y;\eta,c,t)|$ to that of its squarefree restriction. \\ \\
{\bf A reduction to squarefree integers.}
Each $n \in \mathcal{B}_{v}(x,y;\eta,c,t)$ may be decomposed into a powerful and a squarefree component $n= qm$ with $\gcd(q,m)=1$, which is to say that $m$ is squarefree and for each prime $p|q$ we have that $p^2 |q$. Now we may produce an estimate for $|\mathcal{B}_{v}(x,y;\eta,c,t)|$ by studying its squarefree counterpart
\[\mathcal{B}_{\overline{v}}^{\flat}(X,Y;\eta,c,t)=\left\{ m \in S_{\overline{v}}^{\flat} \cap \halfopen{X}{X+Y}: \omega(m_0)<(1-\eta/2)\overline{v}, \ P^{-}(m)>t  \right\}, \] 
where we have written 
\begin{equation}\label{mnought}
    m_0=\prod_{p | m, p \leq (z_{\overline{v} }^{*}(c,X) )^{\eta/3} } p .
\end{equation}
To see why this is true, write $\mathscr{P}(x)$ for the set of powerful numbers not exceeding $x$ and let  $\mathscr{P}_{w}(x)=\mathscr{P}(x) \cap S_{w}(x)$. To deal with the powerful component $q$ of a given integer $n=qm \in \mathcal{B}_{v}(x,y;\eta,c,t)$ let $q_0=\prod_{p^e || q, p \leq z_v^{*}(c,x) } p^e$ and $q_1=q/q_0$ and observe that whenever $q < x_0^{\eta}$ we have $\omega(q_1) < \eta v /2$.  From the assumption that $n=qm \in \mathcal{B}_{v}(x,y;\eta,c,t)$ we get $\omega(q_0)\leq \omega(n_0) \leq (1-\eta)v$ and hence $\omega(q ) \leq  (1- \eta/2)v $. \\
Continuing with the scenario $q < x_0^{\eta}$,  write $X=x/q$,  $Y=y/q$, $w=\omega(q)$ and $w_0=\omega(q_0)$ and consider the decomposition $\overline{v}=\omega(m_0)+\omega(m/m_0)=\overline{v}_0+\overline{v}_1$ with $m_0$ as in \eqref{mnought}. Since $w_0+\overline{v}_0\le (1-\eta)v$ (by definition) and $w_0+\overline{v}_0+\overline{v}_1 \ge (1-\eta/2)v$, we gather that $\overline{v}_1 \ge \eta v/2 \ge \eta \overline{v}/2$ and thus $m \in \mathcal{B}_{\overline{v}}^{\flat}(X,Y;\eta,c,t)$. Employing Golomb's asymptotic \cite{Gol} for powerful integers, we find that 
\begin{align}\label{reductionsfree}
|\mathcal{B}_{v}&(x,y;\eta,c,t)|  \leq \sum_{ \substack{ q m \in \mathcal{B}_{v}(x,y;\eta,c,t)\\ q \in  \mathscr{P}(x_0^{\eta} )   \\ \mu^2(m)=1 }}  1 +  \sum_{q \in  \mathscr{P}(x), q >  x_0^{\eta } }   \big( \frac{y}{q}  +O(1) \big)  \notag \\
&\ll  \sum_{w < (1-\eta/2)v  } \ \sum_{q \in \mathscr{P}_{w}(x_0^{\eta} ) }  \Big| \mathcal{B}_{v -w}^{\flat}( \frac{x}{q},\frac{y}{q} ;\eta,c,t) \Big|  + \min \Big( \frac{y }{ x_0^{\eta/ 2}} , E_{v}(x,y) \Big).
\end{align}
The second remainder follows from the lower bound $\underline{\delta}_{v} (x;t) \gg 1/ (v! \log x)$ \footnote{This lower bound comes from a straightforward modification of the argument for \eqref{lowerrestrictedstar}. } which dominates $x_0^{-\eta/2}$.

For $w \leq (1-\eta/2)v$, we set $\overline{v}=v -w $ and proceed with the treatment of the squarefree bad set. \\

{\bf An estimate for $|\mathcal{B}^{\flat}_{\overline{v}}(X,Y;\eta,c,t)|$.} We first address the range $\exp(2(\log_2 x)^c) \leq v \leq \mathcal{L}_a(x)$. Let $\overline{v} \in [\eta v /2, v]$, let $X \in [x/x_0^{\eta},x] $ and $Y \in [y/x_0^{\eta},y]$. Here we will work with the notation $z_1^{*} =z^{*}_{\overline{v}}(c,X) $ and define
\[\mathcal{I}=[(z_1^{*})^{\eta/3}, z_2^{*} ] , \qquad z_2^{*}= \exp \Big( \frac{ \log X }{\overline{v} \log_3 v} \Big), \qquad r_t= \log \frac{\log z_1^{*} }{\log 2t} \geq 1. \]

Given any $n \in \mathcal{B}^{\flat}_{\overline{v}}(X,Y;\eta,c,t)$ we wish to separate $n$ into three parts $n=n_0 n_1 n_2$ of the form: 
\begin{equation}\label{keyseparation}
n_0 = \prod_{ p|n, \ p < (z_1^{*} )^{\eta/3} } p  , \qquad n_1| \prod_{  p|n, p \in \mathcal{I}  } p=: n_{\mathcal{I}} , \  \qquad n_2=\frac{n}{n_0 n_1}.  
\end{equation} 
Note that any such decomposition guarantees the restriction $n_0 n_1 \leq x^{\varepsilon/4}$. Correspondingly, we set $v_j=\omega(n_j)$ for $j=0,1,2$ so that $\overline{v}= v_0 + v_1+v_2$. The precise choice of $n_1$ will depend on the size of $v_0+\omega(n_{\mathcal{I}} )$. We consider two scenarios:
\begin{itemize}
\item Case I: If $v_0+\omega(n_{\mathcal{I}} )< \overline{v}(1-2/(\log v))$, set $n_1=n_{\mathcal{I}}$.
\item Case II: Otherwise set $k=\lceil \overline{v}(1-2/ \log v) - v_0 \rceil$ and let $n_1$ be the product of the first $k$ prime divisors of $n_{\mathcal{I}}$.
\end{itemize}
Owing to the above restrictions, the summation over $n_2$ will make a small contribution to the cardinality of $\mathcal{B}^{\flat}_{ \overline{v}}(X,Y;\eta,c,t)$. Indeed, in the first case $n_2$ belongs to a set of reduced density, due to the size restrictions on its prime divisors, whereas in the second case we have that $v_2=\omega(n_2) \leq  2 \overline{v} / \log v $. In either scenario we may write  $P^{-}(n_2) > \overline{z}$, where  
\[ \overline{z}=z_2^{*}  \mathds{1}_{ \text{Case I }} + (z_1^{*})^{\eta/3}  \mathds{1}_{ \text{Case II}}. \]

In view of the decomposition \eqref{keyseparation},  we find that
\begin{equation}\label{badbound}
|\mathcal{B}^{\flat}_{ \overline{v}}(X,Y;\eta,c,t)| \leq \sum_{v_0, v_1 }^{'} \  \sum_{ n_0 \in S_{J_0,v_0}(x_0)   } \ 
 \sum_{n_1 \in S_{\mathcal{I}, v_1}(x^{\varepsilon/4 } ) }
  \ \sum_{ \substack{ n_2 \in S_{v_2}(\frac{X}{n_0 n_1},\frac{Y}{n_0 n_1} )  \\ P^{-}(n_2) > \overline{z} }}^{\flat} 1,
\end{equation}
where $J_0=[t, (z_1^{*} )^{\eta/3}] $ and $v_0,v_1$ satisfy the restrictions outlined in Case I and II. Let $m$ denote the product of the first $\lfloor v_2/ 100 \rfloor\ge \overline{v}/ ( 100 \log v  )$ many prime divisors of $n_2$. Since each prime divisor $p |n_2$ exceeds $(z^{*}_1)^{\eta/3}$, this easily implies that $m \asymp U$ for some $U$ in the range \footnote{If necessary,  we may remove some prime divisors of $m$ in order to guarantee the upper bound $m \leq \exp(\log X/ \log_3 X)$.}
\[ \frac{\eta \log X}{ 300 \exp(L_{v}(X)^c ) \log v} \leq \log U \leq \frac{ \log X}{\log_3 X}.\]

Let us write $n_2'=n_2/m$ and set $\omega(n_2')=v_2'$ so that 
\begin{equation}\label{wlower}
 w:=v_0 +v_1+v_2' \geq \overline{v} \big(1- \frac{1}{\log v } \big). 
 \end{equation} 
Before proceeding with the treatment of \eqref{badbound}, we record  the estimate
\begin{align*}
 \Big( \frac{\log X}{v_2' \log \overline{z}} \Big)^{v_2'}  &\leq \Big(\frac{1}{v_2'} \Big)^{v_2'} \Big[  \big( \frac{6}{\eta} v \exp(L_v(x)^c) \big)^{2 \overline{v} / \log v}  +\Big( \frac{\log X}{ \log z_2^{*}  } \Big)^{v_2'} \Big] \\
 &  \leq \Big(\frac{1}{v_2'} \Big)^{v_2'}  \Big[  e^{4 v} +(v \log_3 v )^{v_2'} \Big] 
 \end{align*}  
which is valid since $L_v(x)^c \leq (\log v)/2$.  For fixed $v_0,v_1,v_2'$ and $U$ we thus obtain, via Proposition \ref{extensionprop}, Corollary \ref{Evalcor} and a double application of \eqref{Fnuprescribed},  a contribution no greater than
\begin{align}\label{badsetmainbound}
&    \sum_{ n_0 \in S_{J_0,v_0}(x_0)   } \ 
 \sum_{n_1 \in S_{\mathcal{I}, v_1}(x^{\varepsilon/4 } ) }   \  \sum_{ m \in [U,2U]   }
 \ \ \sum_{  \substack{ n_2' \in S_{v_2'} \\ n_2' \in [\frac{X}{m n_0 n_1},\frac{X+Y}{ m n_0 n_1}]  \\P^{-}(n_2')>\overline{ z} } }  \frac{ F_{v_2'}(n_2') (\log \overline{z})^{-v_2'} }{(v_2') !} \notag \\
& \ll  Y  \frac{ \log x}{(v_2')! } \  \frac{ e^{4 v} \kappa^{v_2'}  + (\kappa v  \log_3 v )^{ v_2'}  }{  (v_2')^{v_2'}}   \sum_{ n_0 \in S_{J_0,v_0}(x_0)   }   \frac{1}{n_0}  
 \sum_{n_1 \in S_{\mathcal{I}, v_1}(x ) }  \frac{1}{n_1} \notag \\
 & \ll Y  \log x  \Big( e^{4 v} + e^{\kappa v /e } \Big)  \frac{ ( \log_3 v )^{ v_2'} }{(v_2')! }  \frac{ (2 L_{v}(x)^c )^{v_1} }{v_1!}
\frac{ (2  r_t )^{v_0}}{v_0 !}
\end{align}
to the RHS of \eqref{badbound}. In the last line we used that $\sup_{v_2' } (\kappa v/v_2')^{v_2'} \leq \exp( \kappa v/e) $. To estimate the final expression in \eqref{badsetmainbound}, we require the inequalities $(a)-(c)$ just below:
\begingroup
\addtolength{\jot}{0.9em}
\begin{align*}
&(a) \qquad   r_{t}^{v_0} L_{v}(x)^{cv_1} (\log_3 v)^{ v_2'}  \leq  (1+o(1))^v r_{t}^{ v - v_1-v_2' } \Big(  L_{v}(x)^c \Big)^{v_1+v_2'  } \\
&(b) \qquad  a_0^{a_0}   a_1^{a_1} a_2^{a_2} \geq (w/3)^w, \qquad \text{ when } a_i \geq 1 \text{ and } a_0+a_1+a_2=w, \\ 
&(c) \qquad w^w \ge ( v/3 )^{v}, \qquad w \in \big[ v(1-\frac{1}{ \log v}), v \big]. 
\end{align*}
\endgroup

\noindent  The three inequalities are left as simple exercises. From the last two bounds combined with a straightforward application of Stirling's approximation, we get that $((v_0 )! (v_1)! (v_2')! )^{-1} \ll \sqrt{ \overline{v}} \ 9^{ \overline{v}} / \overline{v} !$  which, together with the inequality $(a)$ may be inserted into \eqref{badsetmainbound}. Taking into account the number of possible values for $v_0,v_1,v_2'$ and $U$ and invoking \eqref{lowerrestrictedstar} we get that
\begin{align*} 
 |\mathcal{B}^{\flat}_{ \overline{v}}(X,Y;\eta,c,t)|   & \leq \frac{Y (\log x)^5 }{\overline{v} !}  \Big( e^{4 v} + e^{\kappa v /e }  \Big)  9^{ \overline{v} } \big( \frac{3}{2} \big)^v \sum_{w \in [\eta v/3, v]} r_t^{v -w} ( L_w(x)^c)^w  \\
 & \leq Y (\log x)^5  C^{v}  \  \frac{v !}{\overline{v} !}  \  \Delta^{*}_{v,\eta v/3,x'}(x; c,t) , 
  \end{align*}
for some $C =C(\kappa)>0$, any $x'\in [X,X^{3/2}]$ and $\Delta^{*}$ as in \eqref{Deltastardef}. From the first line above we also obtain a direct estimate for the cardinality $ |\mathcal{B}^{\flat}_{ \overline{v}}(X,Y;\eta,c,t)|$. Indeed, since $r_t \leq L_v(x)$, the crude density bound \eqref{deltacrude} reveals that
\begin{equation}\label{secondbadbound}
|\mathcal{B}^{\flat}_{ \overline{v}}(X,Y;\eta,c,t)|   \leq Y (\log x)^6 C^{ v}  \  \frac{v !}{\overline{v} !}  \frac{\delta_v(x)}{L_v(x)^{(1-c) \eta v/3} }.
\end{equation}
 It now follows from \eqref{reductionsfree} that
\begin{align}\label{badboundeta}
 |\mathcal{B}_{\nu}(x,y;\eta,c,t)| & \leq \sum_{w < (1-\eta/2)v  } \ \sum_{q \in \mathscr{P}_{w}(x_0^{\eta }) }  \Big| \mathcal{B}_{v -w}^{\flat}( \frac{x}{q},\frac{y}{q};\eta,c,t) \Big|  + \min \Big( \frac{y }{ x_0^{\eta/ 2}} , E_{v}(x,y) \Big)  \notag \\
& \leq y (\log x)^5    \Delta^{*}_{v, \eta v/3}(x; c,t) C^v  \sum_{w < (1-\eta/2)v  }   \frac{v !}{ (v -w) !}   \ \sum_{q \in \mathscr{P}_{w}(x) } \frac{1}{q}  +E_{v}(x,y) \notag \\
&  \ll  y (\log x)^5   \Delta^{*}_{v,\eta v/ 3}(x; c,t)  C^v \sum_{w < v  } { {v }\choose{ w} }  +E_{v}(x,y) ,
\end{align}
from which we easily deduce the estimate \eqref{badsetstarbound}. Arguing in exactly the same way, the estimate \eqref{secondbadbound} gives rise to \eqref{secondbadsetstarbound}.\\

Finally, when $ \log_3 x \leq v \leq \exp(2 (\log_2 x)^c) $, our task is more straightforward. Set $u=\exp(\varepsilon v_x'/ 4 )$ and separate each $ n \in \mathcal{B}^{\flat}_{ \overline{v}}(X,Y;\eta,c,t)$ as follows:

\begin{equation*}
n_0 = \prod_{ p|n, \ t \leq p \leq (z_1^{*} )^{\eta /3} } p  , \qquad n_1= \prod_{  p|n,  z_1^{*}  \leq p \leq u } p , \  \qquad n_2=\frac{n}{n_0 n_1}.  
\end{equation*} 
Let $v_j=\omega(n_j)$ for $j=0,1,2$. We write $J_0=[t, z_1^{*}]$ and $J_1=[ z_1^{*}, u ] $ and apply \eqref{shortnuvanilla} to find that
\begin{align*}
&|\mathcal{B}^{\flat}_{ \overline{v}}(X,Y;\eta,c,t) | \leq  \sum_{\substack{ v_0,v_1,v_2 \\ v_0 \leq (1-\eta/2) \overline{v} } }  \ \sum_{\substack{  n_0 \in S_{J_0,v_0}(x^{\varepsilon/4})    \\ n_1 \in S_{J_1, v_1}(x^{\varepsilon/4 } )}} \ 
    \ \sum_{  \substack{ n_2 \in S_{v_2} \\ n_2 \in [\frac{X}{ n_0 n_1},\frac{X+Y}{  n_0 n_1}]   } }  \frac{ F_{v_2}(n_2) (\log u)^{-v_2} }{(v_2) !} \\
& \ll  \frac{Y}{ \log x}  \sum_{\substack{ v_0,v_1,v_2 \\ v_0 \leq (1-\eta/2)v  } }   \frac{1}{(v_2 -1)! } \frac{  ( 4v / \varepsilon )^{v_2}}{(v_2)! } \  \sum_{ n_0 \in S_{J_0,v_0}( x^{\varepsilon/4})   }   \frac{1}{n_0}  
 \sum_{n_1 \in S_{J_1, v_1}(x^{\varepsilon/4} ) }  \frac{1}{n_1} 
\end{align*}
and hence, via \eqref{Fnuprescribed} and \eqref{lowerrestrictedstar}  we obtain the estimate
\begin{align*}
& |\mathcal{B}^{\flat}_{ \overline{v}}(X,Y;\eta,c,t) | \ll   \frac{Y}{ \log x} A_{\varepsilon}^v \ \frac{v !}{\overline{v} !}   \sum_{w \in [\eta v/3, v]} r_t^{v -w} ( L_w(x)^c)^w  \\
&  \leq  Y A_{\varepsilon}^v  \  \frac{v !}{\overline{v} !}   \min \Big(  \frac{\delta_v(x)}{L_v(x)^{(1-c) \eta v/2} },  \Delta^{*}_{v, \eta v/3,x'}(x; c,t) \Big)
\end{align*}
for some $A_{\varepsilon} >0$ depending only on $\varepsilon$ and any $x'\in [X,X^{3/2}]$. Inserting this last bound into the first line of \eqref{badboundeta}, we get \eqref{badsetstarbound} and \eqref{secondbadsetstarbound}.
\end{proof}
\end{prop}
\subsection{ Concluding the proof of Theorem \ref{squarefreemainprop}.}

\subsubsection{ The inductive set-up}\label{inductivesetup} Fix $a>4$ and $\varepsilon >0$.  As mentioned previously, we are working under the assumption that $\log_2 x/ (\log_3 x)^{3/4} \leq \nu \leq \log x/(\log_2 x)^{a}$.  By making repeated use of the estimate \eqref{badboundeta}, we aim to prove that almost all integers in $S_{\nu}(x,y)$ are captured by the minorant  $\mathcal{M}_{\nu}^{\sharp} (x,y)$. We will do so by building a hierarchy of collections $\mathcal{G}^{(0)},..., \mathcal{G}^{(t)}$  with the idea that, roughly speaking, each successive $\mathcal{G}^{(r)}$ is made up of integers having fewer large prime divisors than those belonging to its predecessor $\mathcal{G}^{(r-1)}$.  Set
\[ \ell_{\nu}(x) = \frac{ \nu \log(L_{\nu}(x))^2 }{L_{\nu}(x) }, \qquad \eta= ( 10 K_{\varepsilon} ) )^{-2},\qquad  \lambda^{+}=\frac{\log x}{\log_3 x}, \]
where $K_{\varepsilon} >10$ is the constant appearing in Proposition \ref{dealingwiththebadset}. We now define for each index $1 \leq  s \leq \log \nu$ a set of admissible partitions of $\nu$: let $\Gamma^{(s)}_{\nu}$ denote the set of positive integer $(s+2)$-tuples $\underline{w}=(w_0,w_1,...,w_{s+1})$ satisfying the conditions
\begin{itemize}
 \item $w_0 \in \big[ (1-\eta) \nu, \nu \big]$ and  $1 \leq w_j <\eta^{j+1} \nu $ for all $1 \leq j \leq s$ 
 \item $w_j \geq (1-\eta) (\nu-w_0-...-w_{j-1})$ for $1 \leq j \leq s$
 \item $w_0+...+w_{s+1}=  \nu$
 \end{itemize}
We  proceed to factor any given $n \in S_{\nu}(x,y)$ into manageable pieces, keeping in mind the decomposition in \eqref{badsetnudefinition}. For $c\in (0,1)$ and $v \leq \log x$, let us recall the notation
\[  \log z_v^{*}(c,x)=  \frac{ v'_x }{2 \exp( L_v(x)^{c} )}, \qquad \gamma_{v,v_1}(c)= \frac{\rho_v L_{v_1}(x) ^c }{ v_1}, \]
and write 
 \[ d^{*}_{v,c}(m) =\prod_{\substack{ p \ \leq  z_{v}^{*}(c,x)  \\ p^e || m }} p^e , \qquad \  \beta= \beta_{\nu}(x)=\frac{  \log_2(L_{\nu}(x)) }{\log(L_{\nu}(x) ) }.  \]
We will split each $n$ in such a way that 
 \begin{equation} \label{finaldecomp}
 n=n^{(0)}\cdots n^{(s+1)}=n^{-} n^{(s+1)}, \qquad n^{(j)} \in S_{w_j}   
 \end{equation}
 for some admissible $(s+2)$ tuple $\underline{w}=(w_0,...,w_{s+1})$. Starting with the divisor 
 \[ n^{(0)}=d^{*}_{\nu,c_0}(n), \qquad w_0=\omega(n^{(0)}), \qquad c_0=1/2,\] 
 the components $n^{(j)}$ are then defined recursively for all $ 1\leq j \leq s$, setting 
 \[ k_j=\nu-w_0-...-w_{j-1}, \qquad \  n^{(j)}=d^{*}_{k_j,c_j} \Big( \frac{n}{n^{(0)}\cdots n^{(j-1)} }  \Big),  \qquad w_j= \omega(n^{(j)}). \]
Here the choice of each $c_j$ depends on $\nu$ and the preceding coefficients. Having selected $c_0=1/2$, we may define the coefficients belonging to subsequent indices $j \geq 1$  inductively, setting
\begin{equation*} 
c_j=
\left\{
	\begin{array}{lll}
		c_{j-1} \ \ & \mbox{if }  \gamma_{\nu, k_j}(c_{j-1})  \notin [\eta^{2}, \eta^{-2} ]   \\
		 & \\
		c_{j-1} - \beta/2 \ \  & \mbox{otherwise. }   
	\end{array}
\right.
\end{equation*}
This process terminates when after finitely many steps, say $s+1$ steps, we \textbf{either} arrive at a divisor satisfying
 \[ \omega \Big( d^{*}_{k_{s+1},c_{s+1}} \Big( \frac{n}{n^{(0)}\cdots n^{(s)} }  \Big) \Big) \leq (1-\eta) k_{s+1},\]
\textbf{or else} $k_{s+1} \leq \ell_{\nu}(x)$. 
Once the final step has been reached, set $n^{(s+1)}= n/ (n^{(0)}\cdots n^{(s)} )$. In the event that $k_{s+1} > \ell_{\nu}(x)$, we gather that $n^{(s+1)} \in \mathcal{B}_{w_{s+1} } \Big(\frac{x}{n^{-} } , \frac{y}{n^{-}   }; \eta,c_{s+1},z^{*}_{k_s}(c_s)  \Big) $.
\begin{notation}
 In what follows, it will be convenient to discard integers $n\in \halfopen{x}{x+y}$ that have a powerful component $q > \exp(\lambda^{+} /4)$. The integers with this property form an exceptional set $ \mathscr{E}_{\mathscr{P} }= \mathscr{E}_{\mathscr{P} }(x,y)$ of size $O(y/ \exp(-\lambda^{+} /8))$.
\end{notation}
Our next goal is to show that the sequences $\left\{ k_j \right\}, \left\{ c_j \right\}$ and $\left\{  n^{(j)} \right\}$ obtained from the above construction satisfy the following useful properties:
\begin{itemize}
\item The sequence $z^{*}_{k_j}(c_j)$ is non-decreasing for $j=1,...,s$.
\item Setting $k_0=\nu$, we have that $k_{j+1} \leq \eta k_j$ for $j=0,1,..,s$. In particular, $w_{s+1} \leq \eta k_s$.
\item Given any decomposition $n=n^{-} n^{(s+1)}$ as in \eqref{finaldecomp}, the "small component" $n^{-}$ satisfies $P^{+}(n^{-}) \leq z^{*}_{k_s}(c_s)$, whereas $P^{-}(n^{(s+1)} ) > z^{*}_{k_s}(c_s)$. Moreover,  we have that $n^{-} \leq \exp(\lambda^{+})$ for all $n$ \textsl{outside}  $\mathscr{E}_{\mathscr{P} }$.
\item The coefficients $c_j$ take values in $\left\{1/2 -(\beta/2) \N \right\}$. Moreover,  assume that $w_{s+1} \geq \ell_{\nu}(x)$ and write $\gamma=\gamma_{\nu, w_{s+1} }(c_{s+1}) $. Then we have the estimate
\begin{equation}\label{yoverey} \frac{1}{e^{\gamma w_{s+1}/3} } \big[ \big( 9 K_{\varepsilon}  \gamma / \eta \big)^{w_{s+1}} + \big( 9 K_{\varepsilon}  \gamma / \eta \big)^{\eta w_{s+1}/3 }  \big]  \leq 2^{- \eta w_{s+1}}. \end{equation}
\end{itemize}

We will only address the last two assertions since the remaining observations follow immediately from the construction.  In order to establish \eqref{yoverey},  we claim that $\gamma_{\nu,k_j}(c_j) \notin [\eta^2,\eta^{-2}]$ for all $j=1,...,s+1$. To prove the claim, consider two possibilities. Either $\gamma_{\nu,k_j}(c_{j-1}) \notin [\eta^2,\eta^{-2}]$ and hence, by construction, $\gamma_{\nu,k_j}(c_j)=\gamma_{\nu,k_j}(c_{j-1}) \notin [\eta^2,\eta^{-2}]$. Otherwise $\gamma_{\nu,k_j}(c_{j-1}) \in [\eta^2,\eta^{-2}]$ and we have that
\[\eta^2 \leq  \gamma_{\nu,k_j}(c_{j-1})= \frac{\rho_{\nu}(x) L_{k_j}(x) ^{c_{j-1}} }{ k_j}\sim \frac{ \nu L_{k_j}(x) ^{c_{j-1}} }{ k_j L_{\nu}(x) } 
\sim \frac{ \nu L_{\nu}(x) ^{c_{j-1} -1} }{ k_j }.  \]
Since each $k_j \geq \ell_{\nu}(x)$, we infer that $c_{j-1} \geq 3\beta/2$ and thus $c_j \geq \beta$. Moreover, 
\[ \gamma_{\nu,k_j}(c_j)= \gamma_{\nu,k_j}(c_{j-1} -\beta/2) \ll  \eta^{-2}  L_{k_j }(x)^{-\beta/2} <  \eta^2 .  \]
The estimate  \eqref{yoverey} follows immediately from the claim. 

Finally,  we bound the size of $n^{-}$. To treat typical values of $n$, lying outside $\mathscr{E}_{\mathscr{P} }$, we first recall that $k_j \leq 2 \eta^{j+1} \nu$ (by the second bullet point just above) and thus $s \leq \log (L_{\nu}(x))$. Since
\[ L_{k_j}(x)^{c_j}  \geq L_{\nu}(x)^{\beta} \geq  \log (L_{\nu}(x)),\]
it follows that
\begin{align*} 
n^{-}&=n^{(0)} \cdots n^{(s)}  \leq   \exp( \frac{ \lambda^{+} }{4} )  \prod_{j=0}^{s}   z^{*}_{k_j}(c_j)^{w_j}   \leq       \exp \Big(\frac{ \lambda^{+} }{4} +  \log x  \sum_{0 \leq j \leq s} \frac{1}{ \exp(L_{k_j}(x)^{c_j} ) } \Big)   \\ 
&\leq     \exp \Big(\frac{ \lambda^{+} }{4} + \frac{\log x \log (L_{\nu}(x) )}{L_{\nu}(x)} \Big)  \leq \exp(\lambda^{+}). 
\end{align*}

\subsubsection{Assembling the pieces}

For each admissible tuple $\underline{w} \in \Gamma^{(s)}_{\nu}$ we may now collect all the integers of "type $\underline{w}$". Let us separate those tuples $\underline{w}$ for which $w_{s+1} > \ell_{\nu}(x)$ from those satisfying $w_{s+1} \leq \ell_{\nu}(x)$.  Accordingly, we define the sets 
\begin{align*} \mathcal{G}^{(s)}_{\underline{w},+} =\Big \{ n&=n^{(0)}\cdots n^{(s+1)} =n^{-} n^{(s+1)} \in S_{\nu}(x,y) \setminus  \mathscr{E}_{\mathscr{P} }: n^{(j)} \in S_{w_j} \text{ for all $\  1\leq j \leq s$},  \\
& n^{(s+1)} \in \mathcal{B}_{w_{s+1} } \Big(\frac{x}{n^{-} } , \frac{y}{n^{-}   }; \eta,c_{s+1},z^{*}_{k_s}(c_s)  \Big),  \ w_{s+1 }  > \ell_{\nu}(x)  \Big \}
 \end{align*}
and
\begin{equation*} \mathcal{G}^{(s)}_{\underline{w},-} =\Big \{ n=n^{(0)}\cdots n^{(s+1)}  \in S_{\nu}(x,y) \setminus  \mathscr{E}_{\mathscr{P} }: n^{(j)} \in S_{w_j}, \  w_{s+1} \leq \ell_{\nu}(x) \Big \}.
 \end{equation*}
The decomposition  \eqref{finaldecomp} gives rise to a partition (recall the notation in \eqref{goodvsbad} with $c=c_0=1/2$)
\[\mathcal{G}_{\nu}(x,y;\eta,c_0,1) \setminus \mathscr{E}_{\mathscr{P} }= \bigcup_{s} \bigcup_{\underline{w} \in \Gamma^{(s)}_{\nu}} [\mathcal{G}^{(s)}_{\underline{w},+} \cup \mathcal{G}^{(s)}_{\underline{w},-}  ]  \] 
from which we may isolate the collections $\mathcal{G}^{(s)}_{\underline{w},+} $. Indeed, when $w_{s+1} \leq \ell_{\nu}(x) $, any $n=n^{-} n^{(s+1)}  \in \mathcal{G}^{(s)}_{\underline{w},-} $ satisfies $n^{-} \leq \tau$ and thus $n$ must appear in the minorant $\mathcal{M}_{\nu}^{\sharp}(x,y)$.  All that remains is to estimate the sum  
\begin{equation*} 
  I:=  \sum_{s} \sum_{ \ell_{\nu}(x) \leq w \leq \nu} \  \sum_{\substack{    \underline{w} \in \Gamma^{(s)}_{\nu}   \\  w_{s+1} =w }  } |\mathcal{G}^{(s)}_{\underline{w},+} |.
  \end{equation*}
To this end, we will combine \eqref{Gluing}, Proposition \ref{dealingwiththebadset} and \eqref{yoverey}. Writing $\tau=\exp(\lambda^{+})$ and $\ell=\ell_{\nu}(x) $, we first find that
\begin{align*} 
&I  \leq \sum_{c} \sum_{t} \sum_{ \ell \leq w \leq \nu}  \  \sum_{ \substack{ n \in S_{\nu -w }(\tau) \\ P^{+}(n) \leq t }} \big|\mathcal{B}_{w} \Big( \frac{x}{n }, \frac{y}{n}; \eta,c, t  \Big) \big| \\
& \ll     \sum_{c } \sum_{t}  \sum_{ \substack{  \ell \leq w \leq \nu   } }   \Big( y \sum_{ \substack{ n \in S_{\nu -w }(\tau) \\ P^{+}(n) \leq  t }} K_{\varepsilon}^w  \frac{1}{n}  \Delta^{*}_{w, \eta w/3,x} \big( \frac{x}{n}; c,t \big)  +  \sqrt{x/n}+ \frac{y}{n 2^w}  \ \underline{\delta}_{w} (x/n;t)\Big)  \\
& \ll  O(y\delta_{\nu}(x)/2^{\ell} )+ y   \delta_{\nu}(x)   (\log_2x)^3 \sum_{c } \sum_{t}  \sum_{  \ell \leq w \leq \nu   } \Big\{  w^2 \exp( - w \gamma_{\nu,w}/3)  \\
& \times  \Big[ \Big(  9 \gamma_{\nu,w} / \eta  \Big)^{w} +\Big(  9 \gamma_{\nu,w} / \eta  \Big)^{\eta w/3}  \Big] K_{\varepsilon}^w \Big\}=: I^{'} + O(y\delta_{\nu}(x)/2^{\ell} ).
\end{align*}
Here the variables $c,t,w$ satisfy the relations imposed by our construction in section \ref{inductivesetup}. In particular,  $c$ runs over positive values in $\left\{ 1/2 -(\beta/2) \N \right\}$,  $t$ runs over values of the form $z^{*}_k(c)$ with $k \leq \nu$, and $(c,t,w)$ satisfies \eqref{yoverey}.  The estimate $o(y\delta_{\nu}(x) )$ for the remainder term in the second line is obtained via the argument in Lemma \ref{gluglu}: we have that
\[ \sum_{ \substack{ n \in S_{\nu -w }(\tau) \\ P^{+}(n) \leq  t }} \frac{1}{n 2^w} \ \underline{\delta}_{w} (x/n;t) \leq \frac{\delta_{\nu}(x) }{ 2^w}.   \] 
An application of \eqref{yoverey} reveals that
\[  I^{'}  \ll  y   \delta_{\nu}(x)  \nu^4 (\log_2 x)^3  \sum_{  \ell \leq w \leq \nu   }  2^{- \eta w} =O \big(\frac{ y  \delta_{\nu}(x) }{2^{\eta \ell/2}} \big). \]

In conclusion, we have found that $I= O(y \delta_{\nu}(x)/2^{\eta \ell/2})$ and thus
\[|\mathcal{G}_{\nu}(x,y; \eta, c_0,1) \setminus  \mathscr{E}_{\mathcal{P} }|  \leq I +  \mathcal{M}_{\nu}^{\sharp} (x,y)
\leq \mathcal{M}_{\nu}^{\sharp} (x,y) + O \big(\frac{ y  \delta_{\nu}(x) }{2^{\eta \ell/2}} \big). \]
In view of the decomposition \eqref{goodvsbad}, one last application of Proposition \ref{dealingwiththebadset}, together with the lower bound in \eqref{genupperlowerpinuxy},  yields the estimate
\begin{align*}
 \pi_{\nu}(x,y)& \leq | \mathscr{E}_{\mathscr{P} }| +|\mathcal{G}_{\nu}(x,y; \eta, c_0,1) \setminus \mathscr{E}_{\mathscr{P} }| +  |\mathcal{B}_{\nu }(x,y;\eta,c_0,1)| \\
& \leq \mathcal{M}_{\nu}^{\sharp} (x,y)+O \big(\frac{ y  \delta_{\nu}(x) }{2^{\eta \ell/2}} \big) \leq  \mathcal{M}_{\nu}^{\sharp} (x,y)+o( \pi_{\nu}(x,y)), 
\end{align*} 
which completes the proof of Theorem \ref{squarefreemainprop}.

\section{The proofs of Theorem \ref{mainprop} and Corollary \ref{divcor}}\label{finalsection}

\subsection{An asymptotic for  \texorpdfstring{$\pi_{\nu}(x,y)$}{pi{nu}(x,y)} }
As discussed in section \ref{outlinesection},  our goal is to deduce \eqref{shortHR} from the comparative asymptotic $\mathcal{W}^{v,z} (x,y) \sim (y/x)   \mathcal{W}^{v,z} (x)$ proven in Proposition \ref{extensionprop}.  To conclude the argument,  all that remains is to show that the unweighted sums are essentially supported in the sets 
\begin{equation*}
 \mathcal{A}^{\flat}_{v}(x;z)=\left\{n \in S_{v}^{\flat}(x): P^{-}(n)\geq z    \right\},   \qquad \mathcal{A}^{\flat}_{v}(x,y;z)=\mathcal{A}^{\flat}_{v}(x+y;z) \cap \halfopen{x}{x+y}.
\end{equation*}
To this end we require the two comparative estimates for $\delta_{\nu}$ stated at the end of section \ref{deltasection}.
In particular, we will make use of the following special case of \eqref{deltahomothety}:
\begin{equation}\label{densitystability}
\delta_{v}(x) \ll  \delta_{v}(mx), \qquad \text{for} \ \ \  \log m \leq \sqrt{\log x},  \ v \leq (\log x)^{1/3}. 
\end{equation}

\begin{lem}\label{approxdecompglue}
Fix $\varepsilon > 0$.  Let $ x\geq 1$ be sufficiently large, assume that $\log_3 x\leq v \leq  (\log x)^{1/3}/(\log_2 x)^2$ and $(\log x)^4 \leq z \leq \exp( (\log_2 x)^2)$.  Then for each $1 \leq l \leq v$ we have the decomposition
\begin{equation}\label{Wlower}
\sum_{\substack{n_0 \in S_{\nu -l}(\tau) \\ P^{+}(n_0)\leq z }} \mathcal{W}^{ l,z} \big(\frac{x}{n_0},\frac{y}{n_0} \big) =  l! \sum_{\substack{n_0 \in S_{\nu -l}(\tau) \\ P^{+}(n_0)\leq z }}  \big| \mathcal{A}^{\flat}_{l} \big(\frac{x}{n_0},\frac{y}{n_0};z \big) \big| + O\Big( l!  y \frac{\delta_{\nu}(x)  \log x }{\sqrt{z}} \Big) 
\end{equation}
and, similarly, 
\begin{equation}\label{pinuxycomp}
\pi_{\nu}(x,y) \sim  \sum_{l \leq \nu} \ \sum_{\substack{n_0 \in S_{\nu -l}(\tau) \\ P^{+}(n_0)\leq z }}  \big| \mathcal{A}^{\flat}_{l} \big(\frac{x}{n_0}, \frac{y}{n_0};z \big) \big| + O\Big( \frac{y  \delta_{\nu}(x) \log x}{\sqrt{z}} \Big),
\end{equation}
uniformly in the range $x^{17/30 + \varepsilon} \leq y \leq x$.

\begin{proof}
To address \eqref{Wlower}, we may first restrict the inner-most sum $\mathcal{W}^{l,z} (x/n_0,y/n_0) $ to a summation over values having a small powerful component.  Setting $x_0=\exp( \sqrt{\log x})$,  decompose each $n \in  \halfopen{x/n_0}{(x+y)/n_0} $ into a powerful factor $q$ and a squarefree factor $m$. We see that
\begin{equation*} \mathcal{W}^{l,z} (x/n_0,y/n_0) = \sum_{\substack{n=qm \in \halfopen{x/n_0}{(x+y)/n_0} \\ q \in \mathcal{P}(x_0), \ P^{-}(n)>z } } P_l(n) + O(\sqrt{x} + y/\sqrt{x_0}).
\end{equation*}
Next we open the convolution $P_l(n)$ and separate its higher prime-power components (made up of $r$ parts, say) from the remaining prime components. It follows that the LHS of \eqref{Wlower} is 

\begin{align}\label{powersqfreesep}
\begin{split}
&\sum_{\substack{n_0 \in S_{\nu -l}(\tau) \\ P^{+}(n_0)\leq z }} \Bigg[   l!  \big| \mathcal{A}^{\flat}_{l} \big(\frac{x}{n_0},\frac{y}{n_0};z \big) \big| +O\Big (\sum_{r \leq l} \binom{l}{r}  \\
&  \sum_{\substack{n_1,...,n_r \in \mathscr{P}(x) \\ m=n_1 \cdots n_r \leq x_0 \\ \\ P^{-}(n_j)>z} } (l-r)! \  \mathcal{A}^{\flat}_{l-r} \big(\frac{x}{mn_0},\frac{y}{mn_0};z \big)  \Big) \Bigg]=:M_l(x,y) + E_{v,l}(x,y).
\end{split}
\end{align}  
Before moving on to the treatment of $E_{v,l}(x,y)$ we record, by way of the separation argument given just above, a crude bound for long sums.  Indeed, when dealing with long intervals, we may insert into \eqref{powersqfreesep} the estimates \eqref{densitystability} and \eqref{densitycomparisonv} to obtain the bound
\begin{align}\label{crudeWha}
 \sum_{\substack{n_0 \in S_{a}(\tau) \\ P^{+}(n_0)\leq z }} & \mathcal{W}^{ h,z} \big(\frac{x}{n_0} \big) \leq  h! x \ \delta_{a+h}(x) + O\Big ( \sum_{r \leq h} \binom{h}{r} 
  \sum_{\substack{n_1,...,n_r \in \mathscr{P}(x) \\ m=n_1 \cdots n_r \leq x_0 \\ P^{-}(n_j)>z } } (h-r)! \ \frac{x}{m}  \delta_{a+h-r}(\frac{x}{m})  \Big) \notag \\
& \ll x \sum_{0 \leq r \leq h} \binom{h}{r}  (h-r)! \ \delta_{a+h-r}(x)
  \sum_{\substack{n_1,...,n_r \in \mathscr{P}(x) \\ m=n_1 \cdots n_r \leq x_0 \\ P^{-}(n_j)>z } }  1/m \notag\\
& \ll x  h! \delta_{a+h}(x) \sum_{0 \leq r \leq h} \binom{h}{r} \Big( \frac{\log x}{\sqrt{z}} \Big)^r   \ll x  h! \ \delta_{a+h}(x),
\end{align}
for any $a,h\geq 1$ satisfying $a+h \leq v$. 

To bound the error $E_{v,l}(x,y)$, we apply part b of Proposition \ref{extensionprop} together with \eqref{crudeWha} to find that

\begin{align*}
|E_{v,l}(x,y) | & \ll \frac{y}{x} \sum_{r \leq l} \binom{l}{r}    \sum_{\substack{n_0 \in S_{\nu -l}(\tau) \\ P^{+}(n_0)\leq z }} \ \sum_{\substack{n_1,...,n_r \in \mathscr{P}(x) \\ m=n_1 \cdots n_r \leq x_0 \\ \\ P^{-}(n_j)>z} } \mathcal{W}^{ l-r,z} \big(\frac{x}{mn_0} \big) + O_{\varepsilon} \Big(\frac{y}{(\log x )^v} \Big) \\
& \ll  y  \sum_{r \leq l} \binom{l}{r} (l-r)! \sum_{\substack{n_1,...,n_r \in \mathscr{P}(x) \\ m=n_1 \cdots n_r \leq x_0 \\ \\ P^{-}(n_j)>z} } \frac{1}{m}  \delta_{\nu -r}(x/m) + O_{\varepsilon} \Big(\frac{y}{(\log x )^v} \Big) \\
& \ll \delta_{\nu}(x) y  \sum_{r \leq l} \binom{l}{r} (l-r)! \Big( \frac{\log x}{\sqrt{z}} \Big)^r  \ll  l! y \frac{\delta_{\nu}(x)  \log x }{\sqrt{z}} . 
\end{align*}
This concludes the treatment of \eqref{Wlower}.

To prove the approximate decomposition \eqref{pinuxycomp}, we apply Theorem \ref{squarefreemainprop} (to truncate the size of the variable $n_0$) and follow the argument leading up to \eqref{powersqfreesep},  yielding

\begin{align*}
&\pi_{\nu}(x,y) \sim  \sum_{l \leq \nu} \ \sum_{\substack{n_0 \in S_{\nu -l}(\tau) \\ P^{+}(n_0)\leq z }}  \big| \mathcal{A}_{l} \big(\frac{x}{n_0},\frac{y}{n_0};z \big) \big| =  \sum_{l \leq \nu}  \sum_{\substack{n_0 \in S_{\nu -l}(\tau) \\ P^{+}(n_0)\leq z }} \Bigg[  \big| \mathcal{A}^{\flat}_{l} \big(\frac{x}{n_0},\frac{y}{n_0};z \big) \big|  \\
& +O\Big (\sum_{r \leq l} \binom{l}{r}  \sum_{\substack{n_1,...,n_r \in \mathscr{P}(x) \\ m=n_1 \cdots n_r \leq x_0 \\ \\ P^{-}(n_j)>z} } \  \mathcal{A}^{\flat}_{l-r} \big(\frac{x}{mn_0},\frac{y}{mn_0};z \big)  \Big) \Bigg].
\end{align*}  
The remainder term is estimated in the same way as $E_{v,l}(x,y)$ and makes a contribution no greater than $O( y \log x\  \delta_{\nu}(x) / \sqrt{z} )$, as desired.

\end{proof}
\end{lem}

\subsubsection{The proof of \eqref{shortHR}} Let $\nu \leq (\log x)^{1/3}/(\log_2 x)^2$ and set $z=(\log x)^4$.  Combining the approximate decompositions \eqref{pinuxycomp} and \eqref{Wlower} of Lemma \ref{approxdecompglue} and part b of Proposition \ref{extensionprop},  we get that
\begin{align*}
&\pi_{\nu}(x,y) \sim \sum_{l \leq \nu} \frac{1}{l!} \sum_{\substack{n_0 \in S_{\nu -l}(\tau) \\ P^{+}(n_0)\leq z }} \mathcal{W}^{ l,z} \big(\frac{x}{n_0},\frac{y}{n_0} \big) + O\Big( \frac{y  \delta_{\nu}(x) \log x}{\sqrt{z}} \Big) \\
&=\frac{y}{x} \sum_{l \leq \nu} \frac{1}{l!} \sum_{\substack{n_0 \in S_{\nu -l}(\tau) \\ P^{+}(n_0)\leq z }} \mathcal{W}^{ l,z} \big(\frac{x}{n_0} \big) + O\Big( \frac{y  \delta_{\nu}(x) \log x}{\sqrt{z}} \Big) \sim \frac{y}{x} \pi_{\nu}(x).
\end{align*} 

\subsection{A uniform upper bound for \texorpdfstring{$\pi_{\nu}(x,y)$}{pi{nu}(x,y)} } Fix $a>4$ and $\varepsilon>0$.  Let us assume that $(\log_2 x)^{4} \leq \nu \leq \mathcal{L}_a(x)$, since we have already proven \eqref{shortHR} for smaller values of $\nu$, and set $t=\exp( \log x/ (\ell_{\nu}(x) \log_3 x) )$.  In order to bound $\pi_{\nu}(x,y)$ from above, we proceed in two steps.

First,  suppose that we are given parameters $X \in [x/\tau, x]$, $Y \in [y/\tau, y]$, and $\ell_{\nu}(x)/ \log_2 x \leq w \leq \ell_{\nu}(x)$. Further, set $x_0=\exp(\log x/(\log_2 x)^3 )$ and observe that, outside an exceptional set $E \subset \mathcal{A}_w(X,Y;t)$ of cardinality $|E|\ll Y/\sqrt{x_0} \leq Y/(\log x)^{\nu}$,  the powerful component $q$ of any $n \in \mathcal{A}_w(X,Y;t)$ cannot exceed $x_0$ and thus $s=\omega(q) \leq w/ \log_2 x$. \\
Moreover,  for each $n \in \mathcal{A}_w(X,Y;t) \setminus E$, there is a dyadic power $U$ in the range $\log x/ \log _2 x \leq \log U \leq \log x/ \log_3 x$ and a divisor $n_0|n$ satisfying $n_0 \sim U$ and $w_0:=\omega(n_0) \leq w/ \log_2 x$.  Given such a factorisation of $n$, we write $n=n_0 q n'$ with $w'=\omega(n')=w-s-w_0$.  Combining the estimate
\[ \Big( \kappa\frac{\log x}{w' \log t} \Big)^{w'} \leq \Big( \kappa\frac{\ell_{\nu}(x) \log_3 x}{w} \Big)^{w} \leq (\kappa \log_3 x)^{\ell_{\nu}(x)} \leq (1+o(1))^{\nu} \]

with \eqref{curlyFmean}, we have that

\begin{align}\label{genupperprep}
\sum_{\substack{ n  \in S_w( X,Y) \\  P^{-}(n)> t } } 1 & \leq \sum_{\substack{s, w_0 \leq w/ \log_2 x \\ U \text{dyadic} }} \  \sum_{q \in \mathscr{P}_s(x_0) } \sum_{n_0 \sim U} \ \sum_{ \substack{m \in S_{w'}(X/qn_0,Y/qn_0) \\ P^{-}(m) > t } }^{\flat}  \frac{ F_{w'}(n') }{(w')! (\log t)^{w'} } \notag \\
& \ll \frac{Y}{w!} (\log x)^2 \sum_{s, w_0 \leq w }  \Big( \kappa\frac{\log x}{w' \log t} \Big)^{w'}  \sum_{q \in \mathscr{P}(x_0) } \frac{1}{q} \ll Y (\log x)^5 (1+o(1))^{\nu}  \underline{\delta}_{w} (X;t).
\end{align}
Here we have used the straightforward lower bound $\underline{\delta}_{w} (X;t) \gg 1/ (w! \log X)$ and the inequality $(c)$, just below \eqref{badsetmainbound}, to replace $(w')!$ with $w!$.

In a second step, we apply Theorem \ref{squarefreemainprop} to find that

\begin{align*}
\pi_{\nu}(x,y)& \ll \sum_{ \frac{\ell_{\nu}(x) }{\log_2 x} \leq w \leq \ell_{\nu}(x) } \ 
\sum_{ \substack{m \in S_{\nu -w}(\tau) \\ P^{+}(m)\leq t } } \ 
\sum_{\substack{ n  \in S_w( x/m,y/m) \\  P^{-}(n)> t } } 1 \\
& + \sum_{ w < \frac{\ell_{\nu}(x) }{\log_2 x} }  \sum_{ \substack{m \in S_{\nu -w}(\tau) \\ P^{+}(m)\leq t } } \ 
\sum_{\substack{ n  \in S_w( x/m,y/m) \\  P^{-}(n)> t } } 1 =: \pi^{(1)} +  \pi^{(2)}.
\end{align*}
The first sum may be estimated by way of \eqref{genupperprep}. Indeed, since $\nu \geq (\log_2 x)^{4}$ we find that
\begin{align*}  \pi^{(1)} &\ll y (\log x)^5  (1+o(1))^{\nu}  \sum_{  w \leq \ell_{\nu}(x) }\sum_{ \substack{m \in S_{\nu -w}(\tau) \\ P^{+}(m)\leq t } }  \frac{1}{m} \underline{\delta}_{w} (x/m;t) \\
& \leq y  (\log x)^5 \ell_{\nu}(x) (1+o(1))^{\nu} \delta_{\nu}(x)  \ll y   \delta_{\nu}(x) (1+o(1))^{\nu}.
\end{align*}
On the other hand, when $w < \ell_{\nu}(x) / \log_2 x $, we use the trivial observation
\[ \sum_{\substack{ n  \in S_w( x/m,y/m) \\  P^{-}(n)> t } } 1 \leq \frac{y}{m} \ll y w! (\log x) \frac{1}{m}\underline{\delta}_{w} (x/m;t). \]
Recalling the assumption that $\nu \geq (\log_2 x)^{4}$,  we get the bound $w! \log x \ll (1+o(1))^{\nu}$ in this case,  and thus easily retrieve the estimate $ \pi^{(2)} \ll y  \delta_{\nu}(x) (1+o(1))^{\nu}$, thereby completing the proof of \eqref{genupperlowerpinuxy}.

\subsection{Short divisor sums}

\noindent In this final section we address the two remaining assertions stated in the introduction, which is to say, the short divisor sum estimates given in Corollary \ref{divcor}. Let us recall the definition of the $k$-fold divisor function:
\[ 
\tau_k(n) = \sum_{ \substack{ n_1, \dots, n_k \ge 1 \\ \: n_1 \cdots n_k = n }} 1.
\]

As an immediate corollary of the upper bound \ref{genupperlowerpinuxy}, we obtain estimates for divisor sums over squarefree integers in short intervals.

\begin{cor}
Fix $\varepsilon >0, \ a \in (4, 9/2)$ and recall the notation $\mathcal{L}_a(x) =\log x/ (\log_2 x)^{a}$. Then for all sufficiently large  $x \geq 10$ we have the estimate
\begin{equation}\label{shortsfreediv}
 \sum_{ \substack{ n \in \halfopen{x}{x+y} \\ \omega(n)\leq \mathcal{L}_a(x) }}^{\flat} \tau_k(n) \ll y  (\log x)^{2 + (1+o(1))k},
\end{equation}
uniformly in $x^{17/30 + \varepsilon} \leq y \leq x$ and $k \geq 2$. Moreover, under the additional assumption that $2 \leq k \leq \mathcal{L}_{\gamma}(x)$ for some $\gamma > a+1 >5$, we have the sharper bound
\begin{equation}\label{sharpshortsfreediv}
 \sum_{ \substack{ n \in \halfopen{x}{x+y} \\ \omega(n)\leq \mathcal{L}_a(x) }}^{\flat} \tau_k(n) \ll_{\varepsilon} y   \exp \Big(  \frac{\gamma +\varepsilon}{\gamma - 1}k L_k(x)  \Big).
\end{equation}
\begin{proof} Combining the estimates \ref{genupperlowerpinuxy} and \eqref{deltacrude} with the pointwise evaluation $\tau_k(n) = k^{ \omega(n)}$ (for squarefree $n$), we have that
\begin{align*}  
 \sum_{ \substack{ n \in \halfopen{x}{x+y} \\ \omega(n)\leq \mathcal{L}_a(x) }}^{\flat} & \tau_k(n)  \leq \log x \sup_{\nu \leq  \mathcal{L}_a(x) } k^{\nu} \pi_{\nu}(x,y) \ll y \log x \sup_{\nu \leq  \mathcal{L}_a(x) } ((1+o(1) )k )^{\nu} \delta_{\nu}(x)   \\
& \ll  \  y \log x     \sup_{\nu \leq  \mathcal{L}_a(x) }\frac{ ( e_{\nu}(x) k L_{\nu -1}(x) )^{\nu} }{ \nu!}    \ll y (\log x)^2  \sup_{\nu \leq  \mathcal{L}_a(x) } \Big( \frac{  e \ e_{\nu}(x) k \log_2 x  }{ \nu} \Big)^{\nu}.
 \end{align*}
Here $e_{\nu}(x)=1+o(1)$ is the implicit function appearing in \eqref{deltacrude}. We may then recover \eqref{shortsfreediv} from the general inequality $\sup_{x \geq 1} (c/x)^x \leq \exp (c/e)$.

Under the additional assumption that $k \leq \mathcal{L}_{\gamma}(x)$ for some $\gamma > a+1$,  it will be enough to insert (for any small but fixed $\epsilon$) the inequality
\begin{equation}\label{extragammaassump}
\sup_{\nu \leq  \mathcal{L}_a(x) }\frac{ ( e_{\nu}(x) k L_{\nu -1}(x) )^{\nu} }{ \nu!}   \leq  \sup_{\nu \geq  \frac{k }{10 \log_2 x} +1 } \frac{ ( (\frac{\gamma}{\gamma -1} + \epsilon)k L_{k} (x) )^{\nu} }{ \nu!}  +\exp(k)  
\end{equation}
into the preceding argument, to get \eqref{sharpshortsfreediv}. To see why \eqref{extragammaassump} holds, we consider two ranges for $k$. 

{\bf Range I: $2\leq k \leq \exp(\sqrt{\log_2 x})$.} In this case we have that $L_k(x)\sim \log_2 x$ and thus 
\[e_{\nu}(x) k L_{\nu -1}(x) \leq (1+o(1) ) k \log_2 x \sim k L_k(x), \]
for all $\nu \leq  \mathcal{L}_a(x)$, which readily yields \eqref{extragammaassump}. \\

{\bf Range II: $\exp(\sqrt{\log_2 x}) \leq k \leq \mathcal{L}_{\gamma}(x)$.} First, when $\nu \leq k/(4 \log_2 x)$ we see that $( e_{\nu}(x) k L_{\nu -1}(x) )^{\nu} \leq \exp( 2 \nu \log[ k L_{\nu -1}(x)] ) \leq e^k$ so that \eqref{extragammaassump} is satisfied in this particular range of $\nu$.

When $ k/(4 \log_2 x) < \nu \leq  \mathcal{L}_a(x) $, first observe that 
\[ \log_2 \nu \geq \log (\log (k+1) -\log_3 x +O(1))=\log_2(k+1) + o(\log_3 x).  \] 
It follows that $\log( (\nu -1) \log \nu) \geq \log(k \log(k+1)) - (1+o(1)) \log_3 x$  and thus $L_{\nu -1}(x) \leq L_k(x) + (1+o(1)) \log_3 x$. Finally, we see that
\[ \frac{ L_{\nu -1}(x)}{L_k(x)} \leq 1+ (1+o(1)) \frac{\log_3 x }{L_k(x)}  \leq \frac{\gamma +o(1)}{\gamma -1} \]
since the second expression just above, when viewed as a continuous function of $k$, is maximal when $k=\mathcal{L}_{\gamma}(x)$.  This concludes the treatment of \eqref{extragammaassump}.
\end{proof}
\end{cor}

\noindent We are almost ready to establish the short-interval divisor sum bounds \eqref{shortdiv} and \eqref{sharpshortdiv}. The final ingredient in the proof is an estimate for the \textsl{square harmonic} divisor sum \cite[Lemma 2.2]{BN}. 
\begin{lem}
For all $k \geq 2$ and some absolute constant $B>0$, we have that
\begin{equation}\label{harmonicdivisorsquares}
\sum_{m \geq 1} \frac{ \tau_{k}(m^2 )}{m^2 }  \ll (\log(Bk))^{11k}.
\end{equation}

\begin{proof}
The LHS of \eqref{harmonicdivisorsquares} is given by the Euler product 
\begin{equation} 
\sum_{m \geq 1} \frac{ \tau_{k}(m^2 )}{m^2 }=\prod_{p} \left(1+ \sum_{j \geq 1} \frac{\tau_k(p^{2j}) }{p^{2 j}} \right)=:\prod_{p} A_k(p).
\end{equation}
The last displayed sequence of estimates in the proof of \cite[Lemma 2.2]{BN} give the bound
\[ \prod_{p} A_k(p) \leq \exp(Bk) \prod_{p < 7k} [12 k \exp(11k/p)], \]
which readily yields the estimate \eqref{harmonicdivisorsquares} after an application of Mertens' Theorem.
\end{proof}   
\end{lem}

{\bf Proof of Corollary \ref{divcor}.} Let us first address the estimate \eqref{shortdiv}. Set $H=x^{\varepsilon/2}$  and observe that each natural number $m$ admits a factorisation $m=\widetilde{m}\cdot m'$ where $\widetilde{m}$ is its largest square divisor and $m'$ is squarefree.  In view of this factorisation and the pointwise bound 
\[ \tau_k(m) \leq k^{\Omega(m)} \leq   H^{1/4} \  \mathds{1}_{ \text{$k \leq (\log x)^2 $ }} + e^k \ \mathds{1}_{ \text{$k > (\log x)^2 $ }}=:M(k,x), \qquad \ \ \Omega(m)=o\big( \frac{\log x}{ \log_2 x} \big), \] 
together with the sub-multiplicativity of $\tau_k$, we may separate the sum on the LHS of \eqref{shortdiv} to get that   
\begin{align*}
\sum_{\substack{ m \in \halfopen{x}{x+y} \\ \Omega(m)\leq \mathcal{L}_a(x) } }  \tau_{k}(m) \le & 
\sum_{\substack{ m=\widetilde{m} m' \in \halfopen{x}{x+y} \\H \leq \widetilde{m} \leq x } }  M(k,x)
+ \sum_{ \substack{ \widetilde{m} < H \\  \Omega(m)\leq \mathcal{L}_a(x)  }} \tau_k(\widetilde{m}) 
\sum_{\substack{ m' \in  \halfopen{ \frac{x}{\widetilde{m} } }{ \frac{x+y}{\widetilde{m}  }  }   \\ \omega(m') \leq \mathcal{L}_a(x) }}^{\flat} \tau_k(m')\\
&=: M(k,x)  \mathcal{T}_1+\mathcal{T}_2.
\end{align*}
Let us first consider $\mathcal{T}_1$.  Since $x^{17/30}H < y$, we have that
\begin{equation*}
\mathcal{T}_1 \ll\sum_{  \substack{ H \leq \widetilde{m} \leq x  }}^{\square} \Big( 1+ \frac{y}{\widetilde{m}} \Big)   \ll  
 \frac{y}{H^{1/2 }  }+ \sqrt{x} =o\left( \frac{y}{H^{1/4 }} \right).
\end{equation*}
Here, the superscript $\square$ indicates a summation over perfect squares. To deal with $\mathcal{T}_2$, we combine \eqref{harmonicdivisorsquares} and \eqref{shortsfreediv} to find that
\[ \mathcal{T}_2 \ll y (\log x)^{2 + (1+o(1))k} \sum_{  \widetilde{m} < H }^{\square} 
\frac{ \tau_k(\widetilde{m}) }{\widetilde{m}} \ll y  (\log(Bk))^{11k}  (\log x)^{2 + (1+o(1))k}. \]
Collecting the estimates for $\mathcal{T}_1$ and $\mathcal{T}_2$ and applying the uniform inequality $M(k,x) \leq H^{1/4} e^k $, we retrieve \eqref{shortdiv}. 

The proof of \eqref{sharpshortdiv} follows the exact same argument, except that in the treatment of $\mathcal{T}_2$, we use \eqref{sharpshortsfreediv} in place of \eqref{shortsfreediv}.
\qed

\end{document}